\documentclass[a4paper,10pt]{article} \usepackage{amsmath,amssymb,amsthm}
\usepackage[english]{babel}

\newtheorem{theo}{Theorem}[section] \newtheorem{defi}[theo]{Definition}
\newtheorem{lemm}[theo]{Lemma} \newtheorem{prop}[theo]{Proposition}
\newtheorem{coro}[theo]{Corollary}

\newcommand{\Na}{\mathbb N}                   

\newcommand{\Za}{\mathbb Z}                   

\newcommand{\Ra}{\mathbb R}                   

\newcommand{\finpreuve}{\hfill $\Box$}

\newcommand{\name}{$\underline{\qquad \qquad}$}

\newcommand{\refe}[1]{\ref{#1}} \newcommand{\reff}[1]{(\ref{#1})}

\begin{document}

\author{  Jean-Marc
Bouclet\footnote{Jean-Marc.Bouclet@math.univ-lille1.fr}\\
 \\ Universit\'e  de Lille  1 \\ Laboratoire Paul Painlev\'e \\ UMR  CNRS
8524,  \\ 59655 Villeneuve  d'Ascq }
\title{{\bf \sc 
 Littlewood-Paley decompositions on manifolds with ends}}

\maketitle

\begin{abstract} For certain non compact Riemannian manifolds with ends, we obtain Littlewood-Paley type estimates on (weighted)
  $ L^p $ spaces, using the usual square function defined by a dyadic partition.
\end{abstract}

\section{Main results}
\setcounter{equation}{0}

Let $ ({\mathcal M},g) $ be a Riemannian manifold, $ \Delta_g $ the Laplacian on functions and $dg$ the Riemannian measure. Consider a dyadic partition of unit, namely choose $
\varphi_0 \in C_0^{\infty}(\Ra) $ and $ \varphi \in
C_0^{\infty}(0,+\infty) $ such that
\begin{eqnarray}
  1 =
 \varphi_0 (\lambda) + \sum_{k \geq 0} \varphi (2^{-k} \lambda) , \qquad \lambda \geq 0.
 \label{partitiondyadique}
\end{eqnarray}
The existence of such a partition is standard. In this paper, we are 
basically interested in getting estimates of $ ||u||_{L^p ({\mathcal M},dg)} $
in terms of $ \varphi (-2^{-k} \Delta_g)u $, either through the following square function
\begin{eqnarray}
  S_{- \Delta_g} u (\underline{x}) := \left( |\varphi_0 (-\Delta_g) u
 (\underline{x})|^2 +
 \sum_{k \geq 0} | \varphi (-2^{-k} \Delta_g) u(\underline{x})|^2 \right)^{1/2} , \qquad
 \underline{x} \in {\mathcal
M} , \label{fonctioncarreeintro}
\end{eqnarray}
or, at least, through
$$ \left(
\sum_{k \geq 0} || \varphi (-2^k\Delta_{g} ) u
||_{L^p ({\mathcal M },dg)}^2 \right)^{1/2} , $$
and a certain remainder term. For the latter, we think for instance to estimates of the form
\begin{eqnarray}
||  u ||_{L^p ({\mathcal M },dg)} \lesssim \left(
\sum_{k \geq 0} || \varphi (-2^k\Delta_{g} ) u
||_{L^p ({\mathcal M },dg)}^2 \right)^{1/2} + ||u||_{L^2 ({\mathcal M},dg)} , \label{BGTLP}
\end{eqnarray}
for $ p \geq 2 $. In the best possible cases, we want to obtain the  equivalence of norms
\begin{eqnarray}
 ||S_{-\Delta_g}u||_{L^p ({\mathcal M},dg)} \approx || u ||_{L^p ({\mathcal M},dg)} , \label{equivalenceannoncee}
\end{eqnarray}
which is well known, for $ 1 < p < \infty $, if $ {\mathcal M} = \Ra^n $ and $g$ is the Euclidean metric (see for instance \cite{Stei1,Sogg1,Tayl1}).
 
 Such inequalities are typically of interest to localize at high frequencies the solutions (and the initial data) of partial differential equations involving the Laplacian such as the Schr\"odinger equation  $ i \partial_t u  =\Delta_g u $ or the wave equation $ \partial_t^2 u = \Delta_g u $, using that $ \varphi (-h^2 \Delta_g) $ commutes with  $ \Delta_g $. For instance, estimates of the form (\ref{BGTLP})  have been successfully used in \cite{BGT1} to prove Strichartz estimates for the Schr\"odinger equation on compact manifolds. We point out that   the equivalence (\ref{equivalenceannoncee}) actually holds on compact manifolds, but (\ref{BGTLP}) is sufficient to get Strichartz estimates. Moreover  (\ref{BGTLP}) is rather robust and still holds in many cases where (\ref{equivalenceannoncee}) does not. See for instance \cite{Bohyperbolique} where we use this fact.

Littlewood-Paley inequalities on Riemannian manifolds are subjects of
intensive studies.  There is a vast literature in harmonic analysis studying
continuous analogues of the square function (\ref{fonctioncarreeintro}), the
so-called Littlewood-Paley-Stein functions defined via integrals involving the
Poisson  and heat semigroups \cite{Stei1}. An important point is to prove  
$ L^p \rightarrow L^p $ bounds related to these square functions (see for instance
\cite{Loho1} and \cite{Coul1}). However, as explained above, weaker estimates
of the form (\ref{BGTLP}) are often highly sufficient for applications to
PDEs. Moreover, square functions of the form (\ref{fonctioncarreeintro}) are
particularly convenient in microlocal analysis since they involve compactly
supported functions of the Laplacian, rather than fast decaying ones. To
illustrate heuristically this point, we consider the linear Schr\"odinger
equation $ i \partial_t u = \Delta_g u $: if the initial data is spectrally
localized at frequency $ 2^{k/2} $, ie $ \varphi (-2^{-k} \Delta_g) u(0,.)
=u(0,.) $, there is microlocal finite propagation speed stating that the
microlocal support (or wavefront set) of $ u (t,.) $ is obtained by shifting
the one of $ u (0,.) $  along the geodesic flow at speed $ \approx 2^{k/2}
$. This property, which is very useful in the applications, fails if $ \varphi
$ is not compactly supported (away from $ 0 $) and thus reflects the interest
of (\ref{fonctioncarreeintro}). 

As far as dyadic decompositions associated to non constant coefficients operators are
concerned, we have already mentioned \cite{BGT1}. We also have to quote the papers
\cite{KlRo1} and \cite{Schl1}. In \cite{KlRo1}, the authors
develop a dyadic Littlewood-Paley theory   for tensors on compact surfaces 
with limited regularity (of great interest for nonlinear applications). In
\cite{Schl1},  $ L^p $ equivalence of norms  for  dyadic square
functions (including small frequencies) associated to Schr\"odinger operators
are proved for a restricted range of $p$. See also the recent survey \cite{OlZh1} for
Schr\"odinger operators on $ \Ra^n $.

In the present paper, we shall use the analysis of $ \varphi (-h^2 \Delta_g)
$, $ h \in (0,1] $, performed in \cite{Bo2} to derive Littlewood-Paley
inequalities on manifolds with ends (see Definition \ref{defmanifold}). Rather
surprisingly, we were unable to find in the literature a reference for the
equivalence (\ref{equivalenceannoncee}) in reasonable cases such as
asymptotically conical manifolds;  
the latter is certainly clear to specialists. We shall anyway recover this result from our analysis which we can summarize in a model case as follows. Assume for simplicity that a neighborhood of infinity of $ ({\mathcal M}, g) $ is isometric to $ \left(  (R , \infty) \times S , d r^2  + d \theta^2 / w(r)^2 \right) $, with $ (S,d \theta^2) $ compact manifold and $ w (r) > 0 $ a smooth bounded positive function. For instance $ w (r) = r^{-1} $ corresponds to conical ends, and $ w (r) = e^{-r} $ to hyperbolic ends.  We first show  that by considering the modified measure $ \widetilde{dg} = w(r)^{1-n}dg \approx dr d \theta $ and the associated modified Laplacian $ \widetilde{\Delta}_g = w(r)^{(1-n)/2} \Delta_g w (r)^{(n-1)/2} $, we always have the equivalence of norms
$$  ||S_{-\widetilde{\Delta}_g}u||_{L^p ({\mathcal M},\widetilde{dg})} \approx || u ||_{L^p ({\mathcal M},\widetilde{dg})} , $$
for $ 1 < p < \infty $, the square function $ S_{-\widetilde{\Delta}_g} $ being defined by changing $ \Delta_g $ into $ \widetilde{\Delta}_g $
in (\ref{fonctioncarreeintro}). By giving weighted version of this
equivalence, we recover (\ref{equivalenceannoncee}) when $w^{-1}$ is of
polynomial growth. Nevertheless, we emphasize that (\ref{equivalenceannoncee})
can not hold in general for it implies that $ \varphi (- \Delta_g) $ is
bounded on $ L^p ({\mathcal M},dg) $ which may fail for instance in the
hyperbolic case (see \cite{Tayl0}). Secondly, we prove that more robust
estimates of the form (\ref{BGTLP}) always hold and can be spatially localized (see Theorem \ref{theoremedelocalisation}).



Here are the results.

\begin{defi} \label{defmanifold} The manifold $ ({\mathcal M},g ) $ is called almost asymptotic
if  there exists a compact set $ {\mathcal K} \Subset {\mathcal M}
$, a real number $ R  $, a compact manifold $ S $, a function $ r \in C^{\infty}({\mathcal M},\Ra) $ and a function $
w \in C^{\infty}( \Ra , (0,+\infty) ) $  with the following
properties: 
\begin{enumerate}
\item{ $r$ is a coordinate near $ \overline{{\mathcal M} \setminus {\mathcal K}}  $ and
$$ r (\underline{x}) \rightarrow + \infty, \qquad \underline{x} \rightarrow \infty , $$ }
\item{ there is a diffeomorphism
\begin{eqnarray}
  {\mathcal M} \setminus {\mathcal K}  \rightarrow (R,+\infty) \times S, \label{diffeomorphismeabstrait}
\end{eqnarray}
through which  the metric reads in local coordinates
\begin{eqnarray}
 g & = & G_{\rm unif} \left(
r , \theta , d r ,  w(r)^{-1} d\theta \right)
\label{pourlelaplacien} 
\end{eqnarray}
  with 
$$ G_{\rm unif}(r,\theta,V) := \sum_{1 \leq j,k \leq n} G_{jk}(r,\theta) V_j V_k , \qquad V = (V_1, \ldots,V_n)\in \Ra^n , 
$$
if  $ \theta = (\theta_{1}, \ldots , \theta_{n-1}) $ are local
coordinates on $S$.}
\item{ The symmetric matrix $ (G_{jk}(r,\theta))_{1 \leq j,k \leq n} $ has smooth coefficients
such that, locally uniformly with respect to $ \theta $,
\begin{eqnarray}
 \left| \partial_r^j \partial_{\theta}^{\alpha} G_{jk}(r,\theta) \right| & \lesssim & 1, \qquad \qquad r
> R,  \label{uniformegunif} 
\end{eqnarray}
and is uniformly positive definite in the sense that, locally uniformly in $ \theta $,
\begin{eqnarray}
G_{\rm unif}(r,\theta,V) & \thickapprox & |V|^2 , \qquad r
> R,  \ V \in \Ra^{n}. \label{bornespectre}
\end{eqnarray}}
\item{ The function $ w $ is smooth and satisfies,
for all $ k \in \Na $,
\begin{eqnarray}
 w (r) & \lesssim & 1 ,  \label{sanscusp} \\
w (r) / w (r^{\prime}) & \thickapprox & 1 , \qquad \mbox{if} \ \ |r-r^{\prime}| \leq 1 \label{diag}  \\
\left|  d^k w (r) / dr^k \right| & \lesssim & w (r),
\label{cinfiniborne}
\end{eqnarray}
for  $ r,r^{\prime} \in \Ra $. }
\end{enumerate}
\end{defi}

Typical examples are given by asymptotically conical manifolds for which $ w(r) = r^{-1} $ (near infinity) or asymptotically hyperbolic ones for which $ w (r) = e^{-r} $. We note that (\ref{diag}) is equivalent to the fact that, for
some $ C
> 0 $,
\begin{eqnarray}
 C^{-1} e^{-C|r-r^{\prime}|} \leq \frac{w (r)}{  w (r^{\prime}) }
\leq C e^{C|r-r^{\prime}|} . \label{diagequivalent}
\end{eqnarray}
 In particular, this implies that $ w
(r) \gtrsim e^{-C r} $. 

 We recall that, if $ \theta = (\theta_1 , \ldots , \theta_{n-1}) $ are local coordinates on $ S $ and $ (r,\theta) $ are the corresponding ones on $ {\mathcal M} \setminus {\mathcal K} $, the Riemannian measure takes the following form near infinity
$$ d g = w (r)^{1-n} b (r,\theta) dr d \theta_1 \ldots d \theta_{n-1} $$
with $ b (r,\theta) $ bounded from above and from below for $ r \gg 1 $, locally uniformly with respect to $ \theta $ (see \cite{Bo2} for more details).
We also define the density
\begin{eqnarray}
 \widetilde{dg} = w (r)^{n-1} d g  \label{referencedensite}
\end{eqnarray}
and the operator
\begin{eqnarray}
 \widetilde{\Delta}_g = w(r)^{(1-n)/2} \Delta_g w(r)^{(n-1)/2}   . \label{referenceoperateur}
\end{eqnarray}
The multiplication by $ w(r)^{(n-1)/2} $ is  a unitary isomorphism
between $ L^2 ({\mathcal M},\widetilde{dg}) $ and $ L^2 ({\mathcal
M},dg) $ so the operators $ \Delta_g $ and $
\widetilde{\Delta}_g $, which are respectively essentially
self-adjoint on $ L^2 ({\mathcal M},dg) $ and $ L^2 ({\mathcal
M},\widetilde{dg}) $, are unitarily equivalent.

Let us denote by $ P $ either $ - \Delta_{g} $ or $ -
\widetilde{\Delta}_{g} $. For $ u \in
C_0^{\infty}({\mathcal M}) $, we define the square
function $ S_P u $ related to the partition of unit (\ref{partitiondyadique}) by
\begin{eqnarray}
 S_P u (\underline{x}) := \left( |\varphi_0 (P) u
 (\underline{x})|^2 +
 \sum_{k \geq 0} | \varphi (2^{-k} P) u(\underline{x})|^2 \right)^{1/2} , \qquad
 \underline{x} \in {\mathcal
M} . \label{squarefunction}
\end{eqnarray}
We will prove the following result (recall that $ \widetilde{dg} $ and $ \widetilde{\Delta}_g $ are defined by (\ref{referencedensite}) and (\ref{referenceoperateur})).
\begin{theo} \label{theotilde} For all $ 1 < p < \infty $, the
following equivalence of norms holds
$$ || u ||_{L^p ({\mathcal M},\widetilde{dg})} \thickapprox
 || S_{-\widetilde{\Delta}_{g}} u ||_{L^p ({\mathcal M},\widetilde{dg})}
 .  $$
\end{theo}
This theorem implies in particular that $ \varphi_0 (-\widetilde{\Delta}_g) $ and $ \varphi (-2^{k} \widetilde{\Delta}_g) $ are bounded on $ L^p ({\mathcal M},\widetilde{dg}) $. For the Laplacian itself, it is known that compactly supported functions of $ \Delta_g $  are in general not bounded on $ L^p ({\mathcal M},dg) $ (see \cite{Tayl1}) so we can not hope to get the same property. We however have the following result.
\begin{theo} \label{theo} For all $ 2 \leq p < \infty $ and all $
M \geq 0 $,
$$ || u ||_{L^p ({\mathcal M}, dg)} \lesssim
 || S_{-\Delta_{g}} u ||_{L^p ({\mathcal
 M},dg)} + || (- \Delta_{g}+i)^{-M} u ||_{L^2 ({\mathcal M},dg)}
 .
  $$
\end{theo}
 Using the fact that, for $ p \geq 2 $, $ || ( \sum_{k} |u_k|^2 )^{1/2} ||_{L^p} \leq
 (\sum_k ||u_k||_{L^p}^2)^{1/2} $, we
 obtain in particular
\begin{coro} For all $ p \in [2,\infty) $,
\begin{eqnarray}
||u||_{L^p ({\mathcal M},\widetilde{dg})} & \lesssim & \left(
\sum_{k \geq 0} || \varphi (-2^k \widetilde{\Delta}_{g} ) u
||_{L^p ({\mathcal M },\widetilde{dg})}^2 \right)^{1/2} + ||
\varphi_0 (- \widetilde{\Delta}_{g} ) u ||_{L^p ({\mathcal M
},\widetilde{dg})} , \label{pasdeSobolev} \\
||u||_{L^p ({\mathcal M},dg)} & \lesssim & \left( \sum_{k \geq 0}
|| \varphi (-2^k \Delta_{g} ) u ||_{L^p ({\mathcal M },dg)}^2
\right)^{1/2} + ||  u ||_{L^2 ({\mathcal M },dg)}.
\label{SobolevL2}
\end{eqnarray}
\end{coro}
Note the two different situations. In (\ref{SobolevL2}), we have an $ L^2 $ remainder which comes essentially from the Sobolev injection
\begin{eqnarray}
 (1-\Delta_g)^{-n/2 - \epsilon} : L^2 ({\mathcal M},dg) \rightarrow L^{\infty}({\mathcal M}) . \label{injectiondeSobolev}
\end{eqnarray}
 On the other hand in  (\ref{pasdeSobolev}), $ C_0^{\infty} $ functions
of $ \widetilde{\Delta}_{g} $ are bounded on $ L^p ({\mathcal
M},\widetilde{dg}) $, for $ 1 < p < \infty $ (see \cite{Bo2}), but we don't
have Sobolev injections (ie we can not replace $ \Delta_g $ and $dg$ by $ \widetilde{\Delta}_g $ and $ \widetilde{dg} $ in (\ref{injectiondeSobolev})) so we cannot replace $ || \varphi_0
(- \widetilde{\Delta}_{g} ) u ||_{L^p ({\mathcal M
},\widetilde{dg})} $ by $ || u ||_{L^2 ({\mathcal
M},\widetilde{dg})} $.

Actually, we have a result which is more general than Theorem \ref{theotilde}. Consider a temperate weight $ W : \Ra \rightarrow (0,+\infty) $, that is a positive function such that, for some $ C,M > 0 $,
\begin{eqnarray}
 W (r^{\prime}) \leq C W (r) (1+|r-r^{\prime}|)^M, \qquad r,r^{\prime} \in \Ra . \label{poidspolynomial} 
\end{eqnarray}
\begin{theo} \label{theoremeapoids} For all $ 1 < p < \infty $, we have the equivalence of norms
\begin{eqnarray}
|| W(r) u ||_{L^p ({\mathcal M},\widetilde{dg})} \thickapprox
 || W(r) S_{-\widetilde{\Delta}_{g}} u ||_{L^p ({\mathcal M},\widetilde{dg})}
 .  \nonumber
\end{eqnarray}
\end{theo}
This is a weighted version of Theorem \ref{theotilde}. Then, using that
\begin{eqnarray}
 L^p ({\mathcal M},dg) = w(r)^{\frac{n-1}{p}} L^p ( {\mathcal M},
\widetilde{dg} ) , \qquad p \in [1,\infty), \label{weightlp}
\end{eqnarray} 
and that products or (real) powers of weight functions are weight functions, we deduce the following result.
\begin{coro} If $w$ is a temperate weight, then for all $ 1 < p < \infty $, we have the equivalence of norms
$$ || W(r) u ||_{L^p ({\mathcal M}, dg)} \thickapprox
 || W(r) S_{-\Delta_{g}} u ||_{L^p ({\mathcal M},dg)}
 .  $$
\end{coro}
Naturally, this result holds with $ W = 1 $ and we obtain (\ref{equivalenceannoncee}) if $ w $ is a temperate weight. In particular, in the case of
asymptotically euclidean manifolds, this provides a justification of Lemma 3.1 of \cite{BoTz2}.

As noted previously, Theorems \ref{theotilde} and \ref{theo} are interesting to localize some PDEs in frequency. In practice, it is often interesting to localize the datas both spatially and spectrally. For the latter, one requires additional knowledge on the spectral cutoffs, typically  commutator estimates. Such estimates are rather straightforward consequences of the analysis of \cite{Bo2} and allow to prove 
 the following localization property. 
\begin{theo} \label{theoremedelocalisation} Let $ \chi \in C_0^{\infty}({\mathcal M}) $. Assume that $ p \in [2,\infty ) $ and that
\begin{eqnarray}
 0 \leq  \frac{n}{2} - \frac{n}{p} \leq 1 . \label{paireadmissible}
\end{eqnarray}
Then
\begin{eqnarray}
 || (1-\chi) u ||_{L^p ({\mathcal M},dg)} \lesssim \left( \sum_{k \geq 0}
|| (1-\chi) \varphi (-2^k \Delta_{g} ) u ||_{L^p ({\mathcal M },dg)}^2
\right)^{1/2} + ||  u ||_{L^2 ({\mathcal M },dg)}. \label{Schrodinger}
\end{eqnarray}
\end{theo}
This theorem could be generalized by considering for instance more general cutoffs, or even differential operators.
We give only this simple version, which will be used in \cite{Bohyperbolique}
to prove Strichartz estimates at infinity using semi-classical methods in the spirit of \cite{BoTz1}.


\section{The Calder\'on-Zygmund Theorem} \label{section2}
\setcounter{equation}{0}
A basic consequence of the usual Calder\'on-Zygmund theorem is
that  pseudo-differential operators of order $ 0 $ are bounded on $ L^p
(\Ra^n) $ for all $ 1 < p < \infty $. The purpose of this section
is to show a similar result for (properly supported)
pseudo-differential operators with symbols of the form $ a_w
(r,\theta,\rho,\eta) = a (r,\theta,\rho,w(r)\eta) $, with $ a \in
S^0 $. Here the $ L^p $ boundedness will be studied with respect
to the measure $ w (r)^{1-n}dr d\theta $. Recall that $ w $ may
not be bounded from below and hence $ a_w $ doesn't belong to $
S^0 $ in general.

Let us set $ \Omega = (R, + \infty) \times \Ra^{n-1} $ equipped
with the measure
$$ d \nu= w (r)^{1-n} d r d \theta,    $$
$ d \theta $ denoting the Lebesgue measure on $ \Ra^{n-1} $. For
convenience (see Appendix  \refe{analyseharmonique}) and with no
loss of generality, we assume that $ R \in \Na $.
 The
following proposition is a version of the Calder\'on-Zygmund
covering lemma adapted to the measure $d\nu$ (and  to the
underlying metric $ dr^2 + w(r)^{-2} d \theta^2 $).
\begin{prop} \label{lemmeCalderonZygmund}  There exists $ C_0 > 0 $ depending only on $n$ and $w$ such that, for all $ \lambda > 0
$ and all $ u \in L^1 (\Omega,d\nu) $, we can find functions $
\tilde{u}  $,  $ (u_j)_{j \in \Na} $  and a sequence of disjoint
measurable subsets $ ( Q_j )_{j \in \Na} $ of $ \Omega $ such that
\begin{eqnarray}
 & u \ = \  \tilde{u} + \sum_j u_j ,  \label{sommeuj}  & \\
& ||\tilde{u}||_{L^{\infty}(\Omega)} \  \leq \  C_0 \lambda  , & \label{borneponctuelle}   \\
& \int u_j d\nu  \ = \  0 , \qquad \emph{supp}(u_j) \ \subset \ Q_j, & \label{moyennenulle} \\
& \sum_j m (Q_j) \ \leq \  C_0 \lambda^{-1}
||u||_{L^1(\Omega,d\nu)}
, & \label{mesureL1} \\
& || \tilde{u} ||_{L^1 (\Omega,d\nu)} +
  \sum_j || u_j ||_{L^1 (\Omega,d\nu)}  \ \leq \  C_0 ||u||_{L^1(\Omega,d\nu)} . & \label{borneL1}
\end{eqnarray}
The family $ (Q_j)_{j \in \Na}$ can be chosen so that the
following hold : there exist sequences $ (r_j)_{j \in \Na} $, $
(\theta_j)_{j \in \Na}  $ and $ (t_j)_{j \in \Na} $
 such that, for all $ D > 1 $, there exists $ C_D = C (D,n,w) > 0 $ and $ (Q_j^*)_{j \in \Na}
 $ such that, for all $ j \in \Na $:
\begin{eqnarray}
Q_j & \subset & Q_j^* , \\
\nu (Q_j^*) & \leq  & C_D  \nu (Q_j) ,  \label{tjunb}
\end{eqnarray}
and, either
\begin{eqnarray}
 Q_j \subset \left\{ |r-r_j| + \frac{|\theta-\theta_j|}{w(r_j)} \leq t_j
 \right\} \subset \left\{ |r-r_j| + \frac{|\theta-\theta_j|}{w(r_j)} \leq D t_j
 \right\} \subset  Q_j^*  \label{tjun}
\end{eqnarray}
with $ t_j \leq 1 $, or
\begin{eqnarray}
 Q_j \subset \left\{ |r-r_j| \leq 1 \ \mbox{and} \ \frac{|\theta-\theta_j|}{w(r_j)}  \leq t_j \right\}
 \subset \left\{ |r-r_j| \leq 2   \ \mbox{and} \
\frac{|\theta-\theta_j|}{w(r_j)} \leq D t_j \right\}  \subset
Q_j^* , \label{tjde}
\end{eqnarray}
with $ t_j > 1 $.
\end{prop}

\noindent {\it Proof.} See Appendix  \ref{analyseharmonique} .
\finpreuve

In the standard form of this result, each $Q_j $ is a  cube and  $Q_j^*  $ is
its double (obtained by doubling the side of $ Q_j$) and (\ref{tjunb}) can be seen as a consequence
of the usual 'doubling property'. Here the doubling property doesn't hold in
general (typically if $ w (r) = e^{-r} $) but we nevertheless get
(\ref{tjunb}) by replacing (\ref{tjun}) when $t_j$ is large by
(\ref{tjde}). This will be sufficient for we shall use this proposition to
consider operators with
properly supported kernels.

\medskip

Consider next a smooth function $ K $ of the form
$$ K (r,\theta,r^{\prime},\theta^{\prime}) = b \left( r,\theta , r-r^{\prime} , \frac{\theta -
\theta^{\prime}}{w(r)} \right) $$ with $ b $ smooth everywhere and satisfying
\begin{eqnarray}
 \left|  \partial_{\hat{\xi}} b (r,\theta,\hat{\xi})  \right| \leq |\hat{\xi}|^{-1-n}  ,
\qquad  (r,\theta) \in \Omega , \ \hat{\xi} \in \Ra^n \setminus \{ 0 \} , \label{deriveeuniformedeb} \\
\mbox{supp}(b) \subset \Omega \times \{ |\hat{\xi}| < 1 \} . \ \ \
\ \ \qquad \qquad \qquad \label{bdiagonal}
\end{eqnarray}
One must think of $b$ as the Fourier transform $ {\mathcal F}_{\xi
\rightarrow \hat{\xi} } a $ of some symbol $ a (r,\theta,\xi) \in
S^0 (\Ra^n \times \Ra^n) $ (more precisely of some approximation  $ (a_{\epsilon})_{\epsilon \in (0,1]}
$ of $a$ in $
S^{-\infty} $  to ensure the smoothness at $
\hat{\xi} = 0 $), cutoff outside a neighborhood of $ \hat{\xi} = 0 $, and of  $ K $  as the kernel of a
pseudo-differential operator.   We  define the operator $ B$ by
\begin{eqnarray}
 (B u)(r,\theta) = \int_{\Omega} K
(r,\theta,r^{\prime},\theta^{\prime}) u
(r^{\prime},\theta^{\prime}) d \nu (r^{\prime},\theta^{\prime}) .
\label{definitionoperateurpseudow}
\end{eqnarray}
 The assumption $ \reff{bdiagonal} $ states that this operator is
properly supported.
\begin{theo} \label{weak1} There
exists $ C $ such that, for all $B$ as above satisfying the
additional condition
\begin{eqnarray}
 ||B||_{L^2 (\Omega,d\nu) \rightarrow L^2 ({\Omega},d\nu)} \leq 1
, \label{hypotheseborneL2}
\end{eqnarray}
 we have: for all $ u \in L^1 ({\Omega},d\nu) $ and all $
\lambda > 0 $
$$ \nu \left( \{ | B u| > \lambda \} \right) \leq C \lambda^{-1}||u||_{L^1 (\Omega,d\nu)} . $$
 \end{theo}

In other words, $ B $ is of weak type $ (1,1) $ relatively to
$d\nu$.

\medskip

Let us recall a well known lemma on singular integrals.
\begin{lemm} \label{lemmecondHormander} There exists a constant $ c_n $ such that,  for all $ t > 0 $,
 for all $ \tilde{K} \in C^{1}(\Ra^{2n}) $ satisfying
\begin{eqnarray}
 | \partial_y \tilde{K} (x,y)| \leq  |x-y|^{-n-1}, \qquad x \ne y , \ \
x,y \in \Ra^n , \label{condHormander}
\end{eqnarray}
 and for all continuous function
$$ Y : \{|x|> 2 t\} \rightarrow \{|y| < t \} , $$
we have
\begin{eqnarray}
 \int_{|x|>2t} |\tilde{K}(x,Y(x))-\tilde{K}(x,0)| d x \leq c_n  . \label{ydex}
 \end{eqnarray}
\end{lemm}

Note that, in the
usual form of this lemma, the function $ Y $ is simply given by $
Y (x) = y $ with $ |y| < t $ independent of $x$.
Of course, if $ \reff{condHormander} $ is replaced by $ |
\partial_y \tilde{K} (x,y)| \leq  C |x-y|^{-n-1} $ one has to replace $
c_n $ by $ c_n C $ in the final estimate.

For completeness, we recall the simple proof.

\noindent  {\it Proof.} By the Taylor formula and $
\reff{condHormander} $ the left hand side of $ \reff{ydex} $ is
bounded by
$$  \int_{|x| > 2 t} t ||x| - t|^{-n-1} dx = \mbox{vol}({\mathbb S}^{n-1})
\int_{2t}^{\infty} t r^{n-1} (r-t)^{-n-1} dr      $$ where the
change of variable $ u = r/t $ shows that the last integral
is finite and independent of $t$. \finpreuve

\bigskip

\noindent {\it Proof of Theorem  \ref{weak1}.} We use the decomposition (\ref{sommeuj}) and set $ v = \sum_j
u_j $. We have
$$  \nu \left( \{ |B u| > \lambda \} \right) \leq  \nu \left( \{ | B \tilde{u}| > \lambda/2 \} \right) +
\nu \left( \{ |B v| > \lambda / 2 \} \right) . $$ Since 
$$ ||B
\tilde{u}||^2_{L^2 ( \Omega , d\nu)} \leq || \tilde{u}||_{L^2
(\Omega,d\nu)}^2 \leq C_0^2 \lambda ||u||_{L^1(\Omega,d\nu)} , $$
the second inequality being due to (\ref{borneponctuelle}) and (\ref{borneL1}), the Tchebychev inequality yields
$$  \nu \left( \{ | B \tilde{u}| > \lambda/2 \} \right) \leq 4 \lambda^{-2} ||B \tilde{u} ||_{L^2(\Omega,d\nu)}^2
\leq 4 C_0^2 \lambda^{-1} ||u||_{L^1 (\Omega,d\nu)} . $$ We now
have to study $v$. We start by studying the contribution of each
function $ u_j $.

For some fixed $ D > 1 $ large enough to be chosen latter
(independently of $j$), we consider first the situation where $
\reff{tjun} $ holds with $ t_j \leq 1 $. By $ \reff{moyennenulle}
$, we have
$$ \int_{Q_j} K (r,\theta,r^{\prime},\theta^{\prime})
u_j (r^{\prime},\theta^{\prime}) d\nu (r^{\prime},\theta^{\prime}) =
\int_{Q_j} \left( K (r,\theta,r^{\prime},\theta^{\prime}) - K
(r,\theta,r_j,\theta_j) \right) u_j (r^{\prime},\theta^{\prime})
d\nu (r^{\prime},\theta^{\prime})
$$
and thus $ \big| \big| B u_j \big| \big|_{L^1 (\Omega \setminus
Q_j^*,d\nu)}
 \leq ||u_j ||_{L^1(\Omega,d\nu)} \sup_{
 Q_j} I_j (r^{\prime},\theta^{\prime})$, with
\begin{eqnarray}
I_j (r^{\prime},\theta^{\prime}) = \int_{ \Omega \setminus Q_j^* }
\left| K (r,\theta,r^{\prime},\theta^{\prime}) - K
(r,\theta,r_j,\theta_j) \right| d \nu (r,\theta) .
\end{eqnarray}
Using the last inclusion in $ \reff{tjun} $, we get
\begin{eqnarray}
I_j (r^{\prime},\theta^{\prime})  & \leq & \int_{|r-r_j|+
\frac{|\theta - \theta_j |}{w(r_j)} > D t_j} \left| K
(r,\theta,r^{\prime},\theta^{\prime}) - K (r,\theta,r_j,\theta_j)
\right| w (r)^{1-n} dr d \theta  . \label{wr}
\end{eqnarray}
Using successively the changes of variables $ r-r_j \mapsto r $, $
\theta - \theta_j \mapsto \theta $ and $ \theta / w (r+r_j)
\mapsto \theta $, the right hand side of $ \reff{wr} $ can be
written as
\begin{eqnarray}
  \int_{ {\mathcal D}_j } \left| K_j (r,\theta, \tilde{r}_j,
\tilde{\theta}_j) - K_j (r,\theta,0,0) \right|  dr d \theta
\label{vraidomaine}
\end{eqnarray}
with
\begin{eqnarray}
K_j (r,\theta,\tilde{r},\tilde{\theta}) & = & b (r+r_j,w(r+r_j)
\theta
+ \theta_j, r-\tilde{r}, \theta - \tilde{\theta} ) \nonumber \\
\tilde{r}_j & = & r^{\prime} - r_j , \qquad \tilde{\theta}_j \ = \
\frac{\theta^{\prime}-\theta_j}{w(r+r_j)}  . \nonumber \\
{\mathcal D_j} & = & \left\{ (r,\theta) \ | \ |r| + \frac{w(r_j +
r )}{w(r_j)} |\theta|  > D t_j \right\} .
\end{eqnarray}
Recalling that we only consider $ (r^{\prime},\theta^{\prime}) \in
Q_j $ and using $ \reff{tjun} $, we have $ | \tilde{r}_j | \leq
t_j \leq 1 $. Thus, by $ \reff{bdiagonal} $, we have $ |r| \leq 2
$ on both supports of $ K_j (r,\theta, \tilde{r}_j,
\tilde{\theta}_j) $ and $ K_j (r,\theta, 0, 0) $.  By $
\reff{diag} $, this implies that $ w (r_j + r) / w (r_j)
\thickapprox 1 $ on these supports and hence  we can find $ C_1
\geq 1 $, depending only on $ n $ and $w$ such that
$$ \reff{vraidomaine}  \leq \int_{ |r| + |\theta|  > \frac{D}{C_1} t_j }
\left| K_j (r,\theta, \tilde{r}_j, \tilde{\theta}_j) - K_j
(r,\theta,0,0) \right|
 d r d \theta . $$
Now, observe that $ \tilde{\theta}_j \leq t_j w (r_j) / w (r+r_j)
$ with $ |r | \leq 2 $ hence, by possibly increasing $C_1 $, we
also have
$$  | \tilde{r_j} | + | \tilde{\theta}_j | \leq C_1 t_j . $$
By choosing $ D > 2  C_1^2 $, Lemma $ \refe{lemmecondHormander} $
shows that, for all $ (r^{\prime},\theta^{\prime}) \in Q_j $,
$$ \int_{ |r| + |\theta|  > \frac{D}{C_1} t_j }
\left| K_j (r,\theta, \tilde{r}_j, \tilde{\theta}_j) - K_j
(r,\theta,0,0) \right|
 d r d \theta \leq c_n , $$
 since (\ref{deriveeuniformedeb}) implies that (\ref{condHormander}) holds.
 This implies that, if we set $ {\mathcal O} = \cup_{k \in \Na} Q_k^*
 $,
\begin{eqnarray}
 \big| \big| B u_j \big| \big|_{L^1 (\Omega \setminus {\mathcal
O},d\nu)}
 \leq c_n ||u_j ||_{L^1(\Omega,d\nu)} ,
\end{eqnarray}
for all $j$ such that $ \reff{tjun} $ holds with $ t_j \leq 1 $.
We will now prove that this is still true if $ \reff{tjde} $ holds
with $t_j \geq 1 $. To prove the latter, observe first that it is
sufficient to find $D$ large enough such that,
for all such $j$'s,
\begin{eqnarray}
 B u_j (r,\theta) = 0 \qquad \mbox{ for all } \ \ (r,\theta)
\notin Q_j^* \  . \label{grandcube}
 \end{eqnarray}
Indeed, if $ (r^{\prime}, \theta^{\prime}) \in Q_j $, we have $
|r-r_j| \leq |r-r^{\prime}| + 1 $ thus, using $ \reff{bdiagonal}
$, either $ K (r,\theta,r^{\prime},\theta^{\prime}) = 0 $ or $
|r-r_j| \leq 2 $. Assume the latter.  Then, by $ \reff{tjde} $, we
must have $ |\theta -  \theta_j | > D t_j w (r_j) $ for $
(r,\theta)$ outside $ Q_j^* $, hence $ |\theta - \theta^{\prime}|
> (D - 1) t_j w (r_j) $. Since $ |r-r_j| \leq 2 $, $ \reff{diag} $ implies that $
w (r_j)
>  w (r) / C_2 $, for some $C_2 > 0 $ depending only on $w$. By choosing  $ D \geq C_2 + 1 $,
we obtain  $ |\theta - \theta^{\prime}|  > w (r) $ and hence $ K
(r,\theta,r^{\prime},\theta^{\prime}) = 0 $ which completes the
proof of $ \reff{grandcube} $.

Now the conclusion of the proof is standard: we have
\begin{eqnarray}
 \nu \left( \{ |B v| > \lambda / 2 \} \right) & \leq &  \nu ({\mathcal O})
 + \nu \left( \{ |B v| > \lambda / 2 \} \cap \Omega \setminus {\mathcal O} \right)
 \nonumber \\
 & \leq & C_0 C_D \lambda^{-1}||u||_{L^1(\Omega,d\nu)} + c_n
 \lambda^{-1} || v ||_{L^1(\Omega,d\nu)} \nonumber , \\
& \leq & C \lambda^{-1} || u ||_{L^1(\Omega,d\nu)} \nonumber ,
 \end{eqnarray}
 using $ \reff{tjunb} $ and $ \reff{mesureL1} $ to estimate $ \nu ({\mathcal O}) $
 and $ \reff{borneL1} $ for $ || v ||_{L^1(\Omega,d\nu)}  $.
This completes the proof. \finpreuve

\bigskip

The boundedness on $ L^p  $ is then a classical consequence of the
Marcinkiewicz interpolation theorem (see for instance
\cite{Stei1,Tayl1}).

\begin{coro} \label{presqueapplication} For all $ p \in (1,2] $, there exists $ C_p $ such that,
for all $B$ of the form $ \reff{definitionoperateurpseudow} $,
with $ b $ satisfying $ \reff{deriveeuniformedeb} $ and $
\reff{bdiagonal} $, such that $ \reff{hypotheseborneL2} $ holds,
we have
$$ ||B  ||_{L^p (\Omega,d\nu) \rightarrow L^p (\Omega,d\nu)} \leq C_p .
 $$
\end{coro}

\section{Pseudo-differential operators}
\setcounter{equation}{0}

In this section, we apply the results of Section \ref{section2} to the  pseudo-differential operators
involved in the expansions of $ \varphi (-h^2 \Delta_g) $ and $ \varphi (-h^2 \widetilde{\Delta}_g) $ (see \cite{Bo2}).
This means that we consider the following situation. Let
$ ( a_k )_{k \in \Na} $ be a bounded sequence in $
S^{-\infty}(\Ra^n \times \Ra^n) $, ie  for all $ j, l \in \Na $,
all $ \alpha,\beta \in \Na^{n-1} $ and all $ m > 0 $ there exists
$ C $ such that for all $ k \geq 0 $
\begin{eqnarray}
\left| \partial_r^{j} \partial_{\theta}^{\alpha} \partial_{\rho}^l
\partial_{\eta}^{\beta} a_k (r,\theta,\rho,\eta) \right| \leq C_{j l \alpha \beta} (1 + \rho^2 +
|\eta|^2)^{-m} . \label{seminormesuniformes}
\end{eqnarray}
Assume that these symbols are supported in $ \Omega
\times \Ra^n $ where, as in the previous section, $ \Omega = (R,+\infty)
\times \Ra^{n-1}  $.  Fix $ \zeta \in C_0^{\infty}(\Ra^n) $  supported in $ r^2 +
|\theta|^2 <1 $ such that $ \zeta \equiv 1 $ near $0$. For all $ M
\geq 0 $, consider the kernel
\begin{eqnarray}
K_{(M)} ( r,\theta,r^{\prime},\theta^{\prime}) = \sum_{k = 0}^M
2^{kn} \hat{a}_k
 \left( r,\theta,2^k (r-r^{\prime}),2^k \frac{\theta-\theta^{\prime}}{w(r)} \right) \zeta (r-r^{\prime},
\theta-\theta^{\prime}) , \label{kernelexplicite}
\end{eqnarray}
where $ \hat{a} $ is the partial Fourier transform of $ a $  ie $
\hat{a}(r,\theta,\hat{\rho},\hat{\eta}) = \int \! \! \int e^{- i
\rho \hat{\rho}- i \eta \cdot \hat{\eta}} a (r,\theta,\rho,\eta) d
\rho d \eta $. We want to study the boundedness of the associate
operator
\begin{eqnarray}
 (B_{(M)} u)(r,\theta) = \int_{\Omega}
K_{(M)}(r,\theta,r^{\prime},\theta^{\prime}) u
(r^{\prime},\theta^{\prime}) d \nu (r^{\prime},\theta^{\prime}) ,
\label{noyausemblable}
\end{eqnarray}
 on possibly weighted $ L^p $ spaces. Throughout this section, we fix a positive function $ W $ defined on $ \Ra $ such that, for some $ C > 0 $,
\begin{eqnarray}
  W (r) \leq C W (r^{\prime}), \qquad  \mbox{for all} \ r,r^{\prime} \in \Ra \ \mbox{such that} \ |r-r^{\prime}| \leq 1. \label{poidspropre}
\end{eqnarray}
Temperate weights satisfy clearly this condition but as well as powers of $ w
$, although $w$ may not be a temperate weight. 

\begin{lemm} \label{differencefacile} Denote by $ \varrho $ the function
$$ \varrho (r,\theta,r^{\prime},\theta^{\prime}) = \zeta (r-r^{\prime},
\theta-\theta^{\prime}) - \zeta  \left( r-r^{\prime},
\frac{\theta-\theta^{\prime}}{w(r)} \right)  $$ and by $ J_k $ the
function
$$  J_k (r,\theta,r^{\prime},\theta^{\prime}) = 2^{kn} \hat{a}_k
 \left( r,\theta,2^k (r-r^{\prime}),2^k \frac{\theta-\theta^{\prime}}{w(r)} \right) \varrho (r,\theta,
 r^{\prime},\theta^{\prime})  . $$
Define the operator $ R_k $ by
$$ R_k u (r,\theta) =  \int_{\Omega} J_k (r,\theta,r^{\prime},\theta^{\prime}) u (r^{\prime},\theta^{\prime})
 d \nu ( r^{\prime} , \theta^{\prime} ) . $$ Then, for all $
p \in [1,\infty] $,
\begin{eqnarray}
 \sum_{k \geq 0} \left| \left| W(r) R_k W(r)^{-1} \right|
\right|_{L^p (\Omega,d \nu) \rightarrow L^p(\Omega,d \nu) } <
\infty . \label{weightedLp}
\end{eqnarray}
\end{lemm}
Note that $ \reff{weightedLp} $ can be written equivalently as
\begin{eqnarray}
 \sum_{k \geq 0} \left| \left| W(r)w(r)^{\frac{1-n}{p}} R_k W(r)^{-1} w(r)^{\frac{n-1}{p}} \right|
\right|_{L^p (\Omega,dr d \theta) \rightarrow L^p(\Omega,dr d
\theta) } < \infty . \label{weightedLp2}
\end{eqnarray}
using the Lebesgue measure $ d r d \theta $ (with the 
convention that $ (n-1)/p = 0 $ if $ p = \infty$).

\bigskip

\noindent {\it Proof.} Let us prove  (\ref{weightedLp2}). For all $ \gamma \in \Ra $, (\ref{diag}) implies that $ W w^{\gamma} $ also satisfies an estimate of the form (\ref{poidspropre}). We may therefore replace $ W w^{(1-n)/p} $ by $ W $ with no loss of generality.
 Then
$$ (W(r) R_k W(r)^{-1} u )(r,\theta) = \int \! \! \int \widetilde{J}_{k}(r,\theta,r^{\prime},\theta^{\prime}) u (r^{\prime}
,\theta^{\prime}) d r^{\prime} d \theta^{\prime} $$ with 
$$
\widetilde{J}_{k} (r,\theta,r^{\prime},\theta^{\prime}) = w
(r)^{1-n} J_k (r,\theta,r^{\prime},\theta^{\prime}) \times \frac{W(r)}{W(r^{\prime})} . $$
 Since $ \zeta \equiv 1 $ near $ 0 $ and $ w
$ is bounded, there exists  $ c > 0 $ such that, 
\begin{eqnarray}
 |r-r^{\prime}| + \frac{|\theta-\theta^{\prime}|}{w(r)} \geq c , \qquad \mbox{ on  the support
of} \ \varrho . \label{minorationsupportintermediaire}
\end{eqnarray}
 Integrating by part  in the integral
defining $ \hat{a}_k$, one sees that, for all $ N \geq 0 $, $
J_k $ takes the following form
$$ (-1)^N 2^{- ( 2N-n) k} \widehat{ \Delta_{\rho,\eta}^N a}_k \left( r,\theta, 2^{k}(r-r^{\prime}), 2^k
\frac{\theta-\theta^{\prime}}{w(r)} \right) \left(
|r-r^{\prime}|^2 + \frac{|\theta-\theta^{\prime}|^2}{w(r)^2}
\right)^{-N} \varrho (r,\theta,r^{\prime},\theta^{\prime}) .
$$
By the uniform estimates in $k$  (\ref{seminormesuniformes}), (\ref{poidspropre}) and (\ref{minorationsupportintermediaire}), this implies that, for all $N$, there exists $ C_N $ such that
$$ | \widetilde{J}_{k}(r,\theta,r^{\prime},\theta^{\prime}) | \leq C_N 2^{-N k} w (r)^{1-n} \left(
1 + |r-r^{\prime}| + \frac{|\theta-\theta^{\prime}|}{w(r)}
\right)^{-N}
$$ for all $ (r,\theta),(r^{\prime},\theta^{\prime}) \in \Omega $
and all $k \in \Na $. The result  follows then  from the usual
Schur Lemma. \finpreuve

\bigskip

By Lemma \ref{differencefacile}, the $ L^p $ boundedness of $ W(r) B_{(M)} W(r)^{-1} $ is thus
equivalent to the one of $ W(r)\widetilde{B}_{(M)}W(r)^{-1} $ with $ \widetilde{B}_{(M)} $ defined similarly to 
(\ref{noyausemblable})  by the kernel
$$  \widetilde{K}_{(M)} ( r,\theta,r^{\prime},\theta^{\prime}) = \sum_{k = 0 }^{ M} 2^{kn} \hat{a}_k
 \left( r,\theta,2^k (r-r^{\prime}),2^k \frac{\theta-\theta^{\prime}}{w(r)} \right) \zeta \left( r-r^{\prime},
\frac{\theta-\theta^{\prime}}{w(r)} \right) . $$ 
We can then write
$$ (W(r) \tilde{B}_{(M)}W(r)^{-1} u ) (r,\theta) = \int_{\Omega}
b_{M,W}
\left(r,\theta,r-r^{\prime},\frac{\theta-\theta^{\prime}}{w(r)}
\right) u (r^{\prime},\theta^{\prime}) d \nu
(r^{\prime},\theta^{\prime}) $$ where $b_{M,W}$ is defined by
$$ b_{M,W} (r,\theta,\hat{\rho},\hat{\eta}) =
\sum_{k \leq M} 2^{k n} \hat{a}_k (r,\theta,2^k \hat{\rho},2^k
\hat{\eta} ) \zeta (\hat{\rho},\hat{\eta}) \times
\frac{W(r)}{W(r-\hat{\rho})} .
$$

To interpret this operator as an operator of the form (\ref{definitionoperateurpseudow}), with a symbol satisfying (\ref{deriveeuniformedeb}), we need $W$ to be smooth. We thus additionally assume that, for all $ k \geq 0 $,
\begin{eqnarray}
|d^k W (r) / dr^k | \lesssim W(r) . \label{derivabilite}
\end{eqnarray}
We shall see further on that, for the final applications, this condition can doesn't restrict the generality  of our purpose. 
\begin{lemm} \label{conditiondeCalderonZygmund} Assume (\ref{poidspropre}) and (\ref{derivabilite}). There exists $ C > 0 $ such
that, for all $ M \geq 0 $,
\begin{eqnarray}
|\partial_{\hat{\rho},\hat{\eta}} b_{M,W}
(r,\theta,\hat{\rho},\hat{\eta})| & \leq & C (|\hat{\rho}| +
|\hat{\eta}|)^{-n-1} . \label{Szerofacile}
\end{eqnarray}
\end{lemm}

\noindent {\it Proof.} It is standard. We recall it for
completeness. Thanks to the cutoff $ \zeta $, it is sufficient to consider the region where $
|\hat{\rho}|+ |\hat{\eta}| < 1 $.  By
(\ref{seminormesuniformes}) $ \hat{a}_k (r,\theta,.,.) $ is bounded in the Schwartz space as $ (r,\theta)  $ and $k$ vary and, by  (\ref{diag}) and 
(\ref{cinfiniborne}), $ W (r)/W(r -
\hat{\rho}) $ is bounded  on the support of $ \zeta
$ together with its derivatives. Thus, for all $ N > 0 $,
\begin{eqnarray}
 \left| \partial_{\hat{\rho},\hat{\eta}} b_{M,W} (r,\theta,\hat{\rho},\hat{\eta}) \right|  & \leq & C_N
 \sum_{k \geq 0} 2^{k(n+1)}(1+2^{k}|\hat{\rho}| + 2^k
 |\hat{\eta}|)^{-N}
 \nonumber \\
 & \leq & C_N \sum_{k \leq k_0} 2^{k (n+1)} + C_N \sum_{k > k_0} 2^{k (n+1)} 2^{(k_0-k)N}
 \thickapprox C_N 2^{k_0(n+1)} \nonumber
\end{eqnarray}
with $ k_0 = k_0 (\hat{\rho},\hat{\eta})$ such that $ 2^{- k_0-1}
\leq |\hat{\rho}| + |\hat{\eta}| < 2^{ - k_0} $. The result
follows. \finpreuve

\bigskip

We next consider the $ L^2 $ boundedness. 

\begin{lemm} \label{borneL2verifiee}  Assume (\ref{poidspropre}), (\ref{derivabilite}) and the existence of $ C > 1 $ such that, for all $ k \geq C $, we have
\begin{eqnarray}
 (r,\theta,\rho,\eta) \in \emph{supp}(a_k) \Rightarrow C^{-1}
\leq | \rho | + |\eta| \leq C  . \label{presqueorthogonalite}
\end{eqnarray}
Then there exists $ C^{\prime} > 0 $ such that, for
all $ M \geq 0 $,
\begin{eqnarray}
 || W(r) B_{(M)} W(r)^{-1}||_{L^2(\Omega,d\nu)\rightarrow L^2(\Omega,d\nu)} & \leq &
 C^{\prime} .  \label{L2facile}
 \end{eqnarray}
\end{lemm}

\noindent {\it Proof.} The uniform boundedness of the family $ (W(r)
B_{(M)} W(r)^{-1} )_{M \geq 0 } $ on $  L^2
(\Omega,d\nu) )$ is equivalent to uniform boundedness, on $
{\mathcal L}(L^2(\Ra^n,drd\theta)) $, of the family of
pseudo-differential operators with kernels
$$ \int \! \! \int e^{-i (r-r^{\prime})\rho - i (\theta-\theta^{\prime})\cdot \eta }
\widetilde{a}_{M,W} (r,r^{\prime},\theta,\theta^{\prime},\rho,\eta) d
\rho d \eta $$ where
$$ \widetilde{a}_{M,W} (r,r^{\prime},\theta,\theta^{\prime},\rho,\eta) = 
 \frac{W(r)w(r)^{\frac{n-1}{2}}}{W(r^{\prime})w(r^{\prime})^{\frac{n-1}{2}}}
 \zeta
(r-r^{\prime},\theta-\theta^{\prime})
 \sum_{k \leq M} a_k (r,\theta,2^{-k}\rho,2^{-k}w(r)\eta ) . $$
The function in front of the sum is smooth and bounded as
well as its derivatives, by  (\ref{diag}), (\ref{cinfiniborne}), (\ref{poidspropre}), (\ref{derivabilite}) and the compact support of $ \zeta $. The result is then  a
consequence of the Calder\'on-Vaillancourt Theorem since
$$| \partial_x^{\alpha} \partial_{\xi}^{\beta} \sum_{C < k \leq M} a_k(x,2^{-k}\xi) | \leq
C_{\alpha \beta}  \qquad M > 0 , \ (x,\xi) = (r,\theta,\rho,\eta) \in \Ra^{2n } .
$$
This follows from the uniform estimates (\ref{seminormesuniformes}) and the fact that the above sum contains a finite
number of terms, independent of $ x,\xi $ and $ M$ since, by (\ref{presqueorthogonalite}),  $ 2^{-k}|\xi| $ belongs to a fixed compact interval $ [a,b] $  of $  \Ra^+ $ (in particular $ |\xi| \gtrsim 1 $) and
$$ 2^{-k}|\xi| \in [a,b] \qquad \Rightarrow \qquad k \in \left[ \ln_2 |\xi| - \ln_2 b , \ln_2 |\xi | - \ln_2 a \right] $$
where the number of integer points in the last interval is bounded. The proof is complete. \finpreuve

\bigskip


We finally get the following result.

\begin{prop} \label{propositionapoids} Let $ (a_k)_{k \in \Na} $ be a family of symbols supported in $ \Omega \times \Ra^n $ satisfying (\ref{seminormesuniformes})  and (\ref{presqueorthogonalite}). Then, for all positive function $W$ satisfying (\ref{poidspropre})
and all $ p \in (1,2] $, there
exists $ C > 0 $ such that, 
$$ || W(r) B_{(M)} W(r)^{-1} ||_{L^p (\Omega,d\nu) \rightarrow L^p(\Omega,d\nu)} \leq C , \qquad M > 0 . $$
\end{prop}

The following lemma shows that we can assume that $W$ also satisfies (\ref{derivabilite}).

\begin{lemm} We can find $ \widetilde{W} $ satisfying (\ref{poidspropre}), (\ref{derivabilite}) and such that, for some $ C > 1 $, 
\begin{eqnarray}
   W (r)/C \leq \widetilde{W}(r) \leq C W (r) . \label{equivalenceavectilde}
\end{eqnarray}
\end{lemm}

\noindent {\it Proof.} Choose a non zero, non negative $ \omega \in C_0^{\infty}(-1,1) $ and set 
$ \widetilde{W}(r) = \int W (r-s) \omega (s) d s $. Since 
$$ (1+|s|)^{-N} / C \leq W (r-s) / W (r) \leq C (1+|s|)^N , $$
we obtain (\ref{equivalenceavectilde}), which implies in turn that (\ref{poidspropre}) holds for $ \widetilde{W} $ since
$$ \frac{\widetilde{W}(r)}{\widetilde{W}(r^{\prime})} = \frac{\widetilde{W}(r)}{W(r)}  \frac{W(r)}{W(r^{\prime})}  \frac{W(r^{\prime})}{\widetilde{W}(r^{\prime})}   $$
if bounded if $ |r-r^{\prime}| \leq 1 $. This implies
$$ |\widetilde{W}^{(k)}(r) | = | \int W (r-s) \omega^{(k)} (s) d s | \lesssim W (r) \lesssim \widetilde{W}(r) $$
which shows that (\ref{derivabilite}) holds for $ \widetilde{W} $.
\finpreuve

\bigskip

\noindent {\it Proof of Proposition \ref{propositionapoids}.} By (\ref{equivalenceavectilde}), the result holds if and only if it holds with $ \widetilde{W} $ instead of $ W $. We may therefore assume that $ W $ satisfies (\ref{derivabilite}).  By Lemma \ref{borneL2verifiee}, the estimate is true with $ p = 2 $. Then, by Lemma \ref{differencefacile}, it is also true for $ \widetilde{B}_{(M)} $ with $ p =2 $ .
By Lemma \ref{conditiondeCalderonZygmund}, we can apply  Corollary \ref{presqueapplication} to obtain the estimate for all $ 1 < p \leq 2 $ with $ \widetilde{B}_{(M)} $ instead of $ B_{(M)} $
and we conclude using again Lemma \ref{differencefacile}. \finpreuve

\section{Proofs of the main results}
\setcounter{equation}{0}
In  this section,  $ P $ and $  d \mu $
denote either $ - \Delta_{g} $ and $ d g $ or $ -
\widetilde{\Delta}_{g} $ and $ \widetilde{dg} $.
 Using the
partition of unit $ \reff{partitiondyadique} $, we define
$$ A_0 = \varphi_0 (P), \qquad A_k = \varphi (2^{-(k-1)} P), \qquad k \geq 1 , $$
so that,  in the strong sense on $ L^2 ({\mathcal
M },d\mu) $, we have
\begin{eqnarray}
  \sum_{k \geq 0} A_k = 1 , \label{partitionoperateurs}
\end{eqnarray}
and  the square function (\ref{squarefunction}) reads
$$ S_P u (\underline{x}) = ( \sum_k |A_k u(\underline{x})|^2 )^{1/2} , \qquad \underline{x} \in
{\mathcal M} . $$

In the next subsections, we will use the following classical result of harmonic analysis.
 Recall first the definition of the
usual Rademacher sequence $ (f_k)_{k \geq 0} $. For $ k = 0$,  $ f_0 $ is the function
  given on $ [0,1) $ by 
  $$ f_0 (t) = \begin{cases} 1  & \mbox{if} \ 0 \leq t \leq 1/2 \\
 -1 & \mbox{if} \ 1/2 < t < 1 \end{cases} , $$ and  then extended on $
\Ra $ as a $1$ periodic function. If $ k \geq 1 $, $ f_k (t) = f
(2^k t) $, for all $ t \in \Ra $. These functions are orthonormal in $
L^2 ([0,1]) $. Given a sequence of complex numbers $ (a_k)_{k \geq
0} $, if we set
$$ F (t) = \sum_{k \geq 0} a_k f_k (t) , $$ then, for all $ 1 < p < \infty
$, the key estimate related to the Rademacher functions is
\begin{eqnarray}
||F||_{L^2([0,1])} = \big( \sum_{k \geq 0} |a_k|^2 \big)^{1/2}
\leq C_p ||F||_{L^p([0,1])} . \label{Stein}
\end{eqnarray}
For the proof see \cite[p. 276]{Stei1}. As an immediate consequence
of \reff{Stein}, we have the following result.
\begin{prop} \label{carreeutile}
Let $ ( D_k )_{k \geq 0} $ be a family of operators from $
C_0^{\infty}({\mathcal M}) $ to $ L^p ({\mathcal M }, d \mu) $,
for some $ 1 < p < \infty $.
Define the associated square function $ S_D u $ by
$$ S_D u (\underline{x}) = \big( \sum_{k \geq 0} \big|
(D_k u)(\underline{x}) \big|^2  \big)^{1/2}, \qquad \underline{x} \in {\mathcal M} . $$
Then we have
\begin{eqnarray}
 \left| \left| S_D u \right| \right|_{L^p({\mathcal M},d\mu)} \leq C_p \sup_{\tilde{k} \in \Na} \sup_{t \in [0,1]} \big| \big| \sum_{k \leq \tilde{k}  } f_k (t) D_k u \big|
\big|_{L^p ({\mathcal M},d\mu)} . \label{Rademachersupplementaire}
\end{eqnarray}
In particular, if
\begin{eqnarray}
\big| \big| \sum_{k \leq \tilde{k}  } f_k (t) D_k u \big|
\big|_{L^p ({\mathcal M},d\mu)} \lesssim  || u ||_{L^p ({\mathcal
M},d\mu)}, \qquad t \in [0,1], \ u \in C_0^{\infty}({\mathcal M}),
\  \tilde{k}  \geq 0 . \nonumber
\end{eqnarray}
then
\begin{eqnarray}
\left| \left| S_D u \right| \right|_{L^p({\mathcal M},d\mu)}
\lesssim ||u||_{L^p({\mathcal M},d\mu)} , \qquad u \in
C_0^{\infty}({\mathcal M}) . \nonumber
\end{eqnarray}
\end{prop}
\subsection{Proof of Theorems \ref{theotilde} and \ref{theoremeapoids}}
In this part $ P = - \widetilde{\Delta}_g $, $ d \mu = \widetilde{dg} $ and $ W $ is a temperate weight.
\begin{prop} \label{corollairecalculsemiclassique} For all $ N \geq 0 $, we can write
$$ A_k = B_k + C_k , $$
with $ B_k $ such that, for all $ 1 < p \leq 2 $,
$$ \big| \big| \sum_{k \leq \tilde{k}  } f_k (t) W(r) B_k u \big|
\big|_{L^p ({\mathcal M},d\mu)} \lesssim  || W(r) u ||_{L^p ({\mathcal
M},d\mu)}, \qquad t \in [0,1], \ u \in C_0^{\infty}({\mathcal M}),
\  \tilde{k}  \geq 0 , $$
and $ C_k $ such that, for all $ 1 < p < \infty $,
$$ \big| \big| W(r) C_k W(r)^{-1}  \big| \big|_{L^p ({\mathcal M},d\mu) \rightarrow L^p ({\mathcal M},d\mu) 
} \lesssim 2^{-Nk} , \qquad
k \geq 0 . $$
\end{prop}

\noindent {\it Proof.} This follows from the semi-classical parametrix of $ \varphi (-h^2 \widetilde{\Delta}_g) $ given in \cite{Bo2}
 and Proposition \ref{propositionapoids}. 
\finpreuve

\bigskip

We only prove Theorem \ref{theoremeapoids} since Theorem \ref{theotilde} corresponds to the special case $W \equiv 1$. This proof is the standard one to establish the equivalence
of norms of $u$ and $S_P u$ for the usual Littlewood-Paley
decomposition on $ \Ra^n $ (see for instance \cite{Sogg1,Stei1,Tayl1}). We recall it for completeness and to emphasize the difference with the proof of Theorem \ref{theo}.

\bigskip

\noindent {\it Proof of Theorem \ref{theoremeapoids}.} Define $ A_k^W = W(r) A_k W (r)^{-1} $. By Proposition \ref{corollairecalculsemiclassique}, we have
\begin{eqnarray}
\big| \big| \sum_{k \leq \tilde{k}} f_k (t) A_k^W u \big| \big|_{L^p
({\mathcal M}, d \mu) } \lesssim ||u||_{L^p({\mathcal M}, d \mu )}
, \qquad \tilde{k} \geq 0 , \ t \in [0,1], \ u \in
C_0^{\infty}({\mathcal M}) , \label{difficil}
\end{eqnarray}
first for $ 1 < p \leq 2 $, and then for all $ 1 < p < \infty $ by taking the adjoint in the above estimate and replacing $W$ by $ W^{-1} $.
By Proposition \ref{carreeutile},
this implies that
\begin{eqnarray}
 ||W(r) S_P u  ||_{L^p ({\mathcal M}, d \mu )} \lesssim || W(r)  u  ||_{L^p ({\mathcal M}, d \mu )} ,
\qquad u \in C_0^{\infty} ({\mathcal M}) , \label{surLP}
\end{eqnarray}
for $ 1 < p < \infty $. On the other hand, since
$ A_{k_1}A_{k_2} \equiv 0 $ if $ | k_1 - k_2 | \geq 2 $, 
we have
\begin{eqnarray}
 \int_{\mathcal M} \overline{u_1} u_2 d \mu = \sum_{k_1 , k_2 \geq 0 \atop |k_1-k_2|
\leq 1 } \int \overline{A_{k_1}u_1} \ A_{k_2} u_2 \ d \mu . \label{decompositionpresquediagonale}
\end{eqnarray}
  By the Cauchy-Schwarz inequality in the sum,
H\"older's inequality in the integral and (\ref{surLP}) with $ W^{-1} $, we obtain
\begin{eqnarray*}
 \left| \int_{\mathcal M} \overline{u_1} u_2 d \mu \right| & \leq & 3
|| W(r) S_P u_1 ||_{L^p ({\mathcal M}, d \mu)} || W(r)^{-1} S_P u_2
||_{L^{p^{\prime}} ({\mathcal M}, d \mu)}  \\ & \lesssim & || W(r) S_P u_1
||_{L^p ({\mathcal M}, d \mu)} || W(r)^{-1} u_2 ||_{L^{p^{\prime}}
({\mathcal M}, d \mu)}
\end{eqnarray*}
for $ 1 < p  < \infty $, $ p^{\prime} $ being its conjugate
exponent. This then yields the lower bound
$$ || W(r) u_1  ||_{L^p ({\mathcal M}, d \mu)}  \lesssim || W(r) S_P u_1 ||_{L^p ({\mathcal M}, d \mu)},
\qquad u_1 \in C_0^{\infty}({\mathcal M}) , $$ which completes the
proof. \finpreuve






\subsection{Proof of Theorem \ref{theo}}

In this part $ P = - \Delta_g $ and $ d \mu = dg $.
\begin{prop} \label{corollairecalculsemiclassiqueeffectif} For all $ N,M \geq 0 $, we can write
$$ A_k = B_k + C_k , $$
with $ B_k $ satisfying, for all $ 1 < p \leq 2 $,
\begin{eqnarray}
 \big| \big| \sum_{k \leq \tilde{k}  } f_k (t) B_k  u \big|
\big|_{L^p ({\mathcal M},d\mu)}  \lesssim   ||  u ||_{L^p ({\mathcal
M},d\mu)}, \qquad t \in [0,1], \ u \in C_0^{\infty}({\mathcal M}),
\  \tilde{k}  \geq 0 , \label{autuliseravecledual} \\
  || B_k (1-\Delta_g )^{M} ||_{L^2({\mathcal M},d \mu) \rightarrow
L^2({\mathcal M},d \mu)}  \lesssim  2^{ M k/2}, \qquad k \geq 0 ,  \label{Sobolevdeveloppement}
\end{eqnarray}
and $ C_k $ satisfying, for all $ 2 \leq p \leq \infty $,
\begin{eqnarray}
 \big| \big| (1-\Delta_g)^M C_k (1-\Delta_g)^M   \big| \big|_{L^2 ({\mathcal M},d\mu) \rightarrow L^2 ({\mathcal M},d\mu) 
} \lesssim 2^{-Nk} , \qquad
k \geq 0 . \label{Sobolevreste}
\end{eqnarray}
\end{prop}

\noindent {\it Proof.} This follows from the semi-classical parametrix of $ \varphi (-h^2\Delta_g) $ given in \cite{Bo2}
 and Proposition \ref{propositionapoids}. 
\finpreuve

\bigskip

\noindent {\it Proof of Theorem \ref{theo}.} Choose first $ M $ large enough such that we have the Sobolev estimate
\begin{eqnarray}
 \left| \left| (1-\Delta_{g}  )^{-M} \right| \right|_{  L^2
({\mathcal M},dg) \rightarrow L^p ({\mathcal M},dg) } < \infty,
\label{Sobolevfort}
\end{eqnarray}
for all $ 2 \leq p \leq \infty $ (see \cite{Bo2}).
 Denote by $ S_B $ the square function 
 $$ S_B u := (\sum_k
|B_k u|^2 )^{1/2} . $$
Using (\ref{decompositionpresquediagonale}) (with the current new $A_k$) in which we split each $ A_k $ into $ B_k + C_k $, we obtain
\begin{eqnarray}
\left| \int_{\mathcal M} \overline{u_1} u_2 d g \right| & \leq & 3
|| S_B u_1 ||_{L^{p^{\prime}} ({\mathcal M}, d g)} || S_B u_2
||_{L^p ({\mathcal M}, d g)} + ||u_1 ||_{L^{p^{\prime}}
 ({\mathcal M},dg)}  \times \sum_{k_1,k_2 \geq 0 \atop |k_1-k_2| \leq 1 }
  \nonumber \\
 & &
   || B_{k_1}^* C_{k_2} u_2||_{L^p({\mathcal M},dg)} +
  || C_{k_1}^* B_{k_2} u_2 ||_{L^p ({\mathcal M},dg)} +
 || C_{k_1}^* C_{k_2} u_2 ||_{L^p ({\mathcal M},dg)} , \label{SommeSob}
\end{eqnarray}
where $ p^{\prime} $ is the conjugate exponent to $ p $. By (\ref{autuliseravecledual}) applied with $ p^{\prime} $ and
Proposition \ref{carreeutile}, we have $ || S_B u_1
||_{L^{p^{\prime}} ({\mathcal M}, d g)}\lesssim || u_1
||_{L^{p^{\prime}}({\mathcal M},dg)} $. Thus  (\ref{Sobolevdeveloppement}),  (\ref{Sobolevreste}) with $ N \geq M + 1 $, (\ref{Sobolevfort}) and (\ref{SommeSob}) yield
$$ \left| \int_{\mathcal M} \overline{u_1} u_2 d g \right| \lesssim
||u_1 ||_{L^{p^{\prime}}({\mathcal M},dg) } \left(  || S_B u_2
||_{L^p ({\mathcal M}, d g)} + \sum_{k \geq 0} 2^{-Mk/2} ||
(1-\Delta_g)^{-M} u_2 ||_{L^2 ({\mathcal M},dg)} \right) ,  $$
showing that
$$ ||u||_{L^p ({\mathcal M},d\mu)} \lesssim || S_B u
||_{L^p ({\mathcal M}, d g)} +  ||
(1-\Delta_g )^{-M} u ||_{L^2 ({\mathcal M},dg)} . $$
To replace $ S_B $ by $ S_{-\Delta_{g}} $, we introduce the square function
$$ S_C u (\underline{x}) = \left( \sum_{k \geq 0} |C_k u (\underline{x})|^2 \right)^2 $$
so that
\begin{eqnarray}
 || S_B u ||_{L^p ({\mathcal M},dg)} &  \leq &  ||S_{-\Delta_{g}} u ||_{L^p
({\mathcal M},dg)} + || S_C u ||_{L^p ({\mathcal M},dg)} , \nonumber \\ & \lesssim &
  ||S_{-\Delta_{g}} u ||_{L^p ({\mathcal M},dg)} + \sum_{k \geq 0} || C_k u
  ||_{L^p ({\mathcal M},dg)} , \\ & \lesssim &
  ||S_{-\Delta_{g}} u ||_{L^p ({\mathcal M},dg)} + || (1-\Delta_{g})^{-M} u
  ||_{L^2 ({\mathcal M},dg)} ,
\end{eqnarray}
using (\ref{Rademachersupplementaire}), (\ref{Sobolevreste}) with $ N > 0 $ and (\ref{Sobolevfort}). 
\finpreuve





\subsection{Proof of Theorem \ref{theoremedelocalisation}}

We recall first a composition formula for properly supported differential operators. Let $ B_1(h) $ and $ B_2(h) $ be properly supported pseudo-differential operators on $ \Ra^n $ defined by the Schwartz kernels
\begin{eqnarray}
 K_j (x,y,h) = (2 \pi h)^{-n} \int e^{\frac{i}{h} (x-y)\cdot \xi} a_j (x,\xi) d \xi \chi_{j}(x-y), \qquad j = 1,2, \label{noyauAj} 
\end{eqnarray}
where $\chi_j \in C_0^{\infty}(\Ra^n) $, $ \chi_j \equiv 1 $ near $ 0 $ and $ a_j $ symbols in a class that will be specified below. We only assume so far that, for fixed $x$, $ a_j (x,.) $ is integrable.  The kernel $ K (x_1,x_3,h) $ of $ B_1(h) B_2(h) $ is 
\begin{eqnarray}
 (2 \pi h)^{-2n} \int \! \! \int \! \! \int e^{\frac{i}{h} (x_1-x_2)\cdot \xi_1 + \frac{i}{h} (x_2-x_3) \cdot \xi_2}
a_1 (x_1,\xi_1) \chi (x_1-x_2) a_2 (x_2,\xi_2) \chi_2 (x_2-x_3) d \xi_1 d \xi_2 d x_2 , \nonumber 
\end{eqnarray}
that is, using the change of variables $ \xi_1 = \xi_2 + \tau $, $ x_2 = x_1 + t $,
\begin{eqnarray}
 K(x_1,x_3,h)  
 =  (2 \pi h)^{-n} \int  e^{\frac{i}{h} (x_1-x_3)\cdot \xi_2 } a (x_1,x_3,\xi_2,h) d \xi_2 , \nonumber
\end{eqnarray}
with
$$ a (x_1,x_3,\xi_2,h) = (2 \pi h)^{-n} \int \! \! \int  e^{-\frac{i}{h} t \cdot \tau } a_1 (x_1,\xi_2 + \tau) \chi_1 (-t) a_2 (x_1+t,\xi_2) \chi_2 (x_1+t-x_3) dt d \tau . $$
 Since $ \chi_1 $ and $ \chi_2 $ are compactly supported, we can clearly choose $ \chi_3 \in C_0^{\infty}(\Ra^n) $ equal to $1$ near $ 0 $ such that, for all $ v \in \Ra^n $,
$$ \chi_1 (-t) \chi_2 (x_1+t-x_3) = \chi_3 (x_1-x_3)\chi_1 (-t) \chi_2 (x_1+t-x_3)  , $$
which shows that  
$$ K (x_1,x_3,h) = K (x_1,x_3,h) \chi_3 (x_1-x_3) . $$
Assume now that the symbols $ a_j (x,\xi) $ are of the form
\begin{eqnarray}
 a_j (x,\xi) = b_j (r,\theta,\rho,w(r)\eta) , \qquad b_j \in S^{- \infty} (\Ra^n \times \Ra^n) , \label{symboleAj}
\end{eqnarray}
with $ x = (r,\theta) $ and 
$ \xi = (\rho,\eta) $. Writing $ t = (t_r , t_{\theta}) $ and $ \tau = (\tau_{\rho},\tau_{\eta}) $, we then have
$$ a_1 (x_1,\xi_2 + \tau)  a_2 (x_1+t,\xi_2) = b_1 (r_1,\theta_1 , \rho_2 + \tau_{\rho}, w (r_1) (\eta_2 + \tau_{\eta} ) ) 
 b_2 (r_1+ t_r,\theta_1 + t_{\theta}, \rho_2 , w (r_1 + t_r) \eta_2  ) 
$$
which is of the form $ b^{\prime} (r_1,\theta_1,t,\tau,\rho_2,w(r_1)\eta_2) $ with
$$ b^{\prime} (r,\theta,t,\tau,\rho,\eta) =  b_1 (r,\theta , \rho + \tau_{\rho}, \eta + w(r) \tau_{\eta}  ) 
 b_2 \left( r+ t_r,\theta + t_{\theta}, \rho , \frac{ w (r + t_r) }{w(r)} \eta  \right) . $$
Setting
$$ b^{\prime \prime} (r_1,\theta_1,r_3,\theta_3,t,\tau,\rho,\eta) = b^{\prime} (r_1,\theta_1,t,\tau,\rho_2,\eta_2) \chi_1 (-t) \chi_2 ((r_1,\theta_1)-t - (r_3,\theta_3)) , $$
we obtain the following result.
\begin{lemm} Let $B_1(h)$, $ B_2(h) $ be pseudo-differential operators with kernels of the form (\ref{noyauAj}) and with symbols of the form (\ref{symboleAj}). Then the kernel of $ B_1(h) B_2(h) $ is of the form
\begin{eqnarray}
 (2 \pi h)^{-n} \int \! \! \int e^{\frac{i}{h}(r_1-r_3)\rho + \frac{i}{h} (\theta_1 - \theta_3) \cdot \eta}
b (r_1,\theta_1,r_3,\theta_3,\rho,w(r_1)\eta,h) d \rho d \eta \chi_3 \left( (r_1,\theta_1) - (r_3,\theta_3) \right) , \label{noyauacontracter}
\end{eqnarray}
with
$$ b (r_1,\theta_1,r_3,\theta_3,\rho,\eta,h) = (2 \pi h)^{-n} \int \! \! \int  e^{-\frac{i}{h} t \cdot \tau } b^{\prime \prime} (r_1,\theta_1,r_3,\theta_3,t,\tau,\rho,\eta) dt d \tau . $$
\end{lemm}
This result is purely algebraic and becomes  of interest once we have estimates on the symbol $b$. This is the purpose of the next lemma.
\begin{lemm} For all multi-index $ \gamma \in \Na^{3n} $ and all $ m > 0 $
$$ | \partial^{\gamma} b (r_1,\theta_1,r_3,\theta_3,\rho,\eta,h) | \leq C_{m , \gamma} (1+|\rho|+|\eta|)^{-m} . $$
The constant $ C_{m,\gamma} $ depends on a finite number of semi-norms of $b_1$ and $ b_2 $ in $ S^{- \infty} $.
\end{lemm}

\noindent {\it Proof.} This follows standardly from the stationary phase theorem which shows more precisely that $b$ has an asymptotic expansion in powers of $h$whose coefficients are symbols of rapid decay with respect to $ (\rho,\eta) $. \finpreuve

\bigskip

From this lemma and the standard calculus of oscillatory integrals we can write (\ref{noyauacontracter})
with a symbol $ (b_1 \#_w b_2) $ independent of $ (r_3,\theta_3)$ namely,
$$   (2 \pi h)^{-n} \int \! \! \int e^{\frac{i}{h}(r_1-r_3)\rho + \frac{i}{h} (\theta_1 - \theta_3) \cdot \eta}
(b_1 \#_w b_2) (r_1,\theta_1,\rho,w(r_1)\eta,h) d \rho d \eta \chi_3 \left( (r_1,\theta_1) - (r_3,\theta_3) \right) , $$
with an asymptotic expansion in $ S^{- \infty}(\Ra^n \times \Ra^n) $ 
$$  (b_1 \#_w b_2) (r,\theta,\rho,\eta,h) \sim \sum_{k \geq 0} h^k (b_1 \#_w b_2)_k(r,\theta,\rho,\eta) . $$
In particular,
\begin{eqnarray}
 (b_1 \#_w b_2)_0 (r,\theta,\rho,\eta) & = &  b^{\prime \prime}(r,\theta,r,\theta,0,0,\rho,\eta) , \nonumber \\
 & = & b_1 (r,\theta,\rho,\eta) b_2 (r,\theta,\rho,\eta) .
\end{eqnarray} 
This leads to the following result.
\begin{prop} \label{commutateurnonexplicite} Let $B_1(h)$, $ B_2(h) $ be pseudo-differential operators with kernels of the form (\ref{noyauAj}) and with symbols of the form (\ref{symboleAj}). Then 
$$ [B_1(h) , B_2(h)] = h B (h) , $$
where $ B(h) $ has a kernel of the form
$$ (2 \pi h)^{-n} \int \! \! \int e^{\frac{i}{h}(r_1-r_3)\rho + \frac{i}{h} (\theta_1 - \theta_3) \cdot \eta}
c (r_1,\theta_1,\rho,w(r_1)\eta,h) d \rho d \eta \chi_3 \left( (r_1,\theta_1) - (r_3,\theta_3) \right) , $$
$ c (.,h) $ having an asymptotic expansion in $ S^{- \infty}(\Ra^n \times \Ra^n) $. In particular $ (c (.,h))_{h \in (0,1]} $ is bounded in $ S^{- \infty}(\Ra^n \times \Ra^n) $.
\end{prop}

\bigskip

We are now in position to prove Theorem \ref{theoremedelocalisation}. The proof is essentially the same as the one of \cite[Prop. 4.5]{BoTz1}
 using, in the present case, properly supported operators.

\bigskip

\noindent {\it Proof of Theorem \ref{theoremedelocalisation}.} By (\ref{SobolevL2}), we have
$$ ||(1-\chi)u ||_{L^p ({\mathcal M},dg)} \lesssim \left( \sum_{k \geq 0}
||  \varphi (-2^k \Delta_{g} ) (1-\chi) u ||_{L^p ({\mathcal M },dg)}^2
\right)^{1/2} + ||u||_{L^2({\mathcal M},dg)} . $$
Choosing $ \varphi_1 \in C_0^{\infty}(0,+\infty) $ such that $ \varphi_1 \varphi = \varphi $, we can write 
$ \varphi (-h^2 \Delta_{g} ) (1-\chi) $ as the sum of the following three terms
\begin{eqnarray}
Q_1 (h) & = & \varphi_1 (-h^2 \Delta_g)  (1-\chi) \varphi (-h^2 \Delta_g) , \nonumber \\
 Q_2 (h) & = & \left[ \varphi (-h^2 \Delta_g) ,  1-\chi \right] \varphi_1 (-h^2 \Delta_g) , \nonumber \\
 Q_3 (h) & = & \left[  \varphi_1 (-h^2 \Delta_g) , \left[ \varphi (-h^2 \Delta_g) ,  1-\chi \right] \right] . \nonumber
\end{eqnarray}
We recall from \cite{Bo2} that, for each $ N > 0 $ and $ \psi \in C_0^{\infty}(\Ra) $ we can write 
$$ \psi (-h^2 \Delta_g) = B_N (\psi,h) + h^N C_N (\psi,h) $$ with, for all $ q \in [2,\infty] $,
 $$ || C_N (\psi,h) ||_{L^2 ({\mathcal M},dg) \rightarrow L^q ({\mathcal M},dg)} \lesssim  h^{N} , $$
and
$$ || B_N (\psi,h) ||_{L^q ({\mathcal M},dg) \rightarrow L^q ({\mathcal M},dg)} \lesssim 1, \qquad 
|| B_N (\psi,h) ||_{L^2 ({\mathcal M},dg) \rightarrow L^q ({\mathcal M},dg)} \lesssim h^{-n (\frac{1}{2}-\frac{1}{q})} . $$
We recall more precisely that  pseudo-differential operators as those considered in Proposition \ref{commutateurnonexplicite} satisfy estimates as above.
Therefore, using Proposition \ref{commutateurnonexplicite} and (\ref{paireadmissible}), we also have
$$ ||\left[  B_N (\varphi,h) ,  1-\chi \right]  ||_{L^2 ({\mathcal M},dg) \rightarrow L^p ({\mathcal M},dg)} \lesssim 1 . $$
Thus, by choosing $ N $ large enough
\begin{eqnarray}
|| Q_1 (h) u ||_{L^p ({\mathcal M},dg)} & \lesssim & || (1-\chi) \varphi (-h^2 \Delta_g) ||_{L^p ({\mathcal M},dg)} + h || u ||_{L^2 ({\mathcal M},dg)} , \nonumber \\
||  Q_2 (h) u ||_{L^p ({\mathcal M},dg)} & \lesssim & || \varphi_1 (-h^2 \Delta_g) u ||_{L^2 ({\mathcal M},dg)} , \nonumber \\
|| Q_3 (h) u  ||_{L^p ({\mathcal M},dg)} & \lesssim & h || u ||_{L^2 ({\mathcal M},dg)} , \nonumber
\end{eqnarray}
applying Proposition \ref{commutateurnonexplicite} twice for the last estimate, namely that
$$  ||\left[  B_N (\varphi_1,h) , \left[  B_N (\varphi,h) ,  1-\chi \right] \right] ||_{L^2 ({\mathcal M},dg) \rightarrow {L^p ({\mathcal M},dg)}} \lesssim h .$$
 The result then follows from
$$ \sum_{ h = 2^{-k}, k \in \Na} || \varphi_1 (-h^2 \Delta_g) u ||_{L^2 ({\mathcal M},dg)}^2 \lesssim ||  u ||^2_{L^2 ({\mathcal M},dg)} $$
by almost orthogonality and the Spectral Theorem. \finpreuve

\appendix

\section{Proof of Proposition $ \bf \refe{lemmeCalderonZygmund} $} \label{analyseharmonique}

We first define special families of partitions of $ \Omega $.
Given $ n_0 \in \Na $ and $ k  \geq -n_0   $ integer, we denote by
$ {\mathcal P}(k) $ a countable partition of $ \Omega $, i.e.
$$ {\mathcal P}(k) = ({\mathcal P}_{l
}(k))_{ l \in \Na  } , \qquad \Omega = \sqcup_{l \in
\Na} {\mathcal P}_{l} (k) . $$
In the sequel, the sets  $ {\mathcal P}_{\gamma}(k) $ will always be measurable.

 Given a family of partitions $ {\mathcal P}:=({\mathcal
P}(k))_{k \geq -n_0} $, we shall say that:
\begin{itemize}
\item{${\mathcal P}$ is non increasing if, for all  $ k \geq 1 - n_0 $ and all $ l \in \Na $, there
exists  $ l^{\prime} \in \Na $ such that \begin{eqnarray}
 {\mathcal
P}_{l} (k) \subset {\mathcal P}_{l^{\prime}} (k-1),
\label{decro}
\end{eqnarray} }
\item{${\mathcal P}$ is locally finite if, for all compact subset $ K \subset \Omega $ there exists a compact
subset $ K^{\prime} \subset \Omega $ such that, for all $ k \geq -
n_0 $,
\begin{eqnarray}
\bigsqcup_{l \in \Na \atop {\mathcal P}_{l}(k) \cap K
\ne \emptyset} {\mathcal P}_{l}(k) \subset K^{\prime}
\label{locafini}
\end{eqnarray}}
\item{ ${\mathcal P}$ is of vanishing diameter if there exists a sequence $ \epsilon_k \rightarrow 0 $
}such that, for all $ k \geq - n_0 $ and all $ l \in \Na $ there
exists $ x_{k, l} \in \Omega $ such that
$$ {\mathcal P}_{l}(k) \subset \{ x \in \Omega \ | \ |x-x_{k,l}| \leq \epsilon_k \} . $$
\end{itemize}

\bigskip

\noindent {\bf Remarks. 1.} If $ {\mathcal P} $ is non increasing,
in  \reff{decro}, $ l^{\prime} $  is uniquely defined by $ l $ and
$ k $.

\noindent {\bf 2.} If $ {\mathcal P} $ is non increasing, then it
is locally finite if and only if for all compact subset $ K $
there exists another compact subset $ K^{\prime} $ such that
\reff{locafini}  holds for $ k = - n_0 $.

\noindent {\bf 3.} \label{remarque} If $ {\mathcal P} $ is non
increasing, it follows by a simple induction that if  $ {\mathcal
P}_{l_1}(k_1) \cap {\mathcal P}_{l_2} (k_2) \ne
\emptyset $ for some $ k_2 \geq k_1 \geq - n_0 $ and $ l_1,
l_2 \in \Na  $, then $ {\mathcal P}_{l_2}(k_2) \subset
{\mathcal P}_{l_1} (k_1) $.

\bigskip

\begin{defi} \label{definitionremarque} The family of partitions $ ({\mathcal P}(k))_{k \geq -n_0}
$ is admissible if it is non increasing, locally
finite and of vanishing diameter.
\end{defi}

The proof of Proposition $ \refe{lemmeCalderonZygmund}  $ is based
on a suitable choice of admissible partitions which we now
describe.

\medskip

\noindent {\bf Construction of a family of  admissible
partitions.} For $ m = (m_1 , \cdots , m_{n - 1}) \in \Za^{n-1} $,
we set $$ {\mathcal C}_m = [ m_1 , m_1 + 1 ) \times \cdots \times
[ m_{n-1}, m_{n-1} + 1 ) $$ and for $ \tau
> 0 $, we note $  \tau {\mathcal C}_m = \{ \tau \theta \ ; \ \theta \in C_m  \} $ so that $
\sqcup_{m \in \Za^{n-1}} \tau {\mathcal C}_m = \Ra^{n-1}   $ is a
decomposition of $ \Ra^{n-1} $ into cubes of side $ \tau $.
Setting  $ k_+ = \max (0,k) $, we can define, for all $ k \in \Za
$
\begin{eqnarray}
{\mathcal P}_{(i,m)} (k) =   2^{-k_+}(i,i+1] \times 2^{-k} w
([2^{-k_+}i]) {\mathcal C}_m \qquad i \in  \Na \cap [ 2^{k_+} R ,
\infty ), \ m \in \Za^{n-1}, \label{exemple2}
\end{eqnarray}
where $ [2^{-k_+}i] $ denotes the integer part of $ 2^{-k_+}i $.

For notational convenience, we then relabel $ ( {\mathcal
P}_{(i,m)}(k))_{(i,m) \in \Na \cap [2^{k_+}R,\infty) \times
\Za^{n-1}} $ as $ ({\mathcal P}_l (k))_{l \in \Na} $.

Let us notice that, for $ k \in \Za $ and $ l \in \Na $, we have
$$
 \nu ({\mathcal P}_{l}(k)) = 2^{-k(n-1)} \int_{2^{-k_+}i}^{2^{-k_+}(i+1)}
\left( \frac{w([2^{-k_+}i])}{w(r)} \right)^{n-1} \ dr .
$$
Thus, using $ \reff{diag} $, there exists $ C_2 \geq 1 $ such that
\begin{eqnarray}
  C_2^{-(n-1)} 2^{-k(n-1)}
\leq \nu ({\mathcal P}_{l}(k)) \leq C_2^{n-1} 2^{-k(n-1)} \qquad
\mbox{if} \ \ k \leq 0, \label{minoration1} \\
C_2^{-(n-1)} 2^{-kn} \leq \nu ({\mathcal P}_{l}(k)) \leq C_2^{n-1}
2^{-kn} \ \ \ \   \ \   \qquad \mbox{if} \ \ k \geq 1 .
\label{minoration2}
\end{eqnarray}

\begin{lemm} For all $ n_0 \in \Na $, $ ({\mathcal P}(k))_{k \geq - n_0}
\equiv \big(
 ({\mathcal P}_l (k))_{l \in \Na} \big)_{k \geq -n_0} $ (defined by \reff{exemple2}) is an
 admissible family of partitions of $ \Omega $. Furthermore, there
 exists $ C_3 > 1 $ independent of $n_0$ such that, for all $ k \geq 1 - n_0 $ and all $ l
 \in \Na $,
\begin{eqnarray}
C_3^{-1} \leq \frac{\nu ({\mathcal P}_{l}(k))}{\nu({\mathcal
P}_{l^{\prime}}(1-k))} \leq C_3,
 \label{constante2bis}
\end{eqnarray}
with $ l^{\prime} = l^{\prime}(k,l) $ the unique integer
satisfying \reff{decro}.
\end{lemm}

Let us already point out that our  family of admissible partitions
has been designed in order to get  \reff{constante2bis} which will
be crucial in the proof of Lemma \refe{fonctionsCalderonZygmund}
below.

\medskip

\noindent {\it Proof.} For each $ k \geq -n_0 $,  $ {\mathcal
P}(k) =  ({\mathcal P}_l (k))_{l \in \Na}   $ is obviously  a
partition of $ \Omega  $. Since $ w $ is bounded, the family $
{\mathcal P} =  ({\mathcal
  P}(k))_{k \geq -n_0}  $ is of vanishing diameter ($ \epsilon_k \approx
2^{-k}  $). Assuming that it is non increasing, it also easy to
check that $ {\mathcal P}  $  is locally finite and hence
admissible. So let us prove  that  $ {\mathcal P}  $  is non
increasing. If $ k \leq 0 $, we have
$$ (i,i+1] \times 2^{-k} w (i) {\mathcal C}_m  \subset
(i^{\prime},i^{\prime}+1] \times 2^{1-k} w (i^{\prime}) {\mathcal
C}_{m^{\prime}} $$ provided $ i = i^{\prime} $ and $ {\mathcal
C}_m \subset 2 {\mathcal C}_{m^{\prime}} $, which clearly holds
for some $ m^{\prime} \in \Za^{n-1} $. Thus  \reff{decro}  holds
if $ k \leq 0  $. If $ k \geq 1 $, we remark that if
\begin{eqnarray}
 2^{-k} (i,i+1 ] \subset 2^{1-k} ( i^{\prime} , i^{\prime} + 1 ]
 \label{doublei}
 \end{eqnarray}
 then $
   [ 2^{-k} i ] = [ 2^{1-k} i^{\prime} ]  $. This follows easily
    from the fact that $ 2^{-k} (i,i+1) \cap \Na =
 2^{1-k}(i^{\prime},i^{\prime}+1) \cap \Na =
\emptyset $. Thus
$$ 2^{-k}(i,i+1] \times 2^{-k} w ([2^{-k}i]) {\mathcal C}_m  \subset
2^{1-k}(i^{\prime}, i^{\prime} + 1] \times 2^{1-k} w
([2^{1-k}i^{\prime}]) {\mathcal C}_{m^{\prime}} $$ with  $
i^{\prime} $ such that  \reff{doublei}  holds and $ m^{\prime} $
such that $ {\mathcal C}_m \subset 2 {\mathcal C}_{m^{\prime}} $.
Therefore $ {\mathcal P} $ is non increasing.

The estimate  \reff{constante2bis}  follows easily from
\reff{minoration1} and  \reff{minoration2}.  \finpreuve

\bigskip

We now recall a basic result which is a version of Lebesgue's
Lemma.

\begin{lemm} \label{standard} Let $ u \in L^1(\Omega, d\nu) $ and $ {\mathcal P}
$ be an admissible family of partitions of $ \Omega $. Assume that
$ A \subset \Omega $ is a measurable subset such that there exists
$ C > 0  $ satisfying:  for all $  k \geq - n_0 $ and all $ l
\in \Na $
\begin{eqnarray}
  {\mathcal P}_{l} (k) \cap A \ne  \emptyset \ \
\Rightarrow \ \  \frac{1}{\nu({\mathcal P}_{l}(k))}
\int_{{\mathcal P}_{l}(k)} |u(x)| \ d\nu(x) \leq C .
\label{majomoye}
\end{eqnarray}
Then $ |u| \leq C $ almost everywhere on $ A $.
\end{lemm}

\noindent {\it Proof.} For all $ v \in L^1 (\Omega,d\nu) $,
we set
$$ ({\mathcal E}_{k} v)(x) = \sum_{l \in \Na} \frac{1}{\nu (
{\mathcal P}_{l}(k) ) } \int_{ {\mathcal P}_{l}(k) }
v(y) d \nu(y)  \chi_{ {\mathcal P}_{l}(k) }(x) , $$ $
\chi_{ {\mathcal P}_{l}(k) } $ being the characteristic
function of $ {\mathcal P}_{l}(k) $. We first remark that
\begin{eqnarray}
\lim_{k \rightarrow \infty} {\mathcal E}_k v = v \qquad \mbox{in}
\ \ L^1 (\Omega, d\nu). \label{densite1}
\end{eqnarray}
Indeed, since  $ || {\mathcal E}_k v ||_{L^1} \leq
||v||_{L^1} $ for all $v$, we may assume that $ v$ is continuous and compactly supported.
Then, denoting by $ K $ the support of $v$, we have for all $ k
\geq - n_0 $
$$ || {\mathcal E}_k v - v ||_{L^1(\Omega, \rm{d}\nu)} \leq \sum_{l \in \Na \atop
{\mathcal P}_{l}(k) \cap K \ne \emptyset} \nu ( {\mathcal
P}_{l}(k) ) \sup_{x,y \in {\mathcal P}_{l}(k) }
|v(y)-v(x)| .
$$
Using the local finiteness of $ {\mathcal P} $ and the fact that
it is of vanishing diameter, there exists a compact subset $
K^{\prime} $ such that
$$ || {\mathcal E}_k v - v ||_{L^1(\Omega, \rm{d}\nu)} \leq \nu (K^{\prime})
\sup_{x,y \in K^{\prime} \atop |x-y| \leq 2 \epsilon_k}|v(y)-v(x)|
\rightarrow 0, \qquad k \rightarrow \infty,
$$
and  \reff{densite1}  follows. In particular, $ \chi_A {\mathcal E
}_k |u| \rightarrow \chi_A |u| $ in $ L^1 (\Omega,d\nu) $ so there
exists a subsequence $ \chi_A {\mathcal E}_{k_j}|u| $ converging
almost everywhere to $ \chi_A |u| $. Using  \reff{majomoye}  we
have
$$ 0 \leq \left( \chi_A {\mathcal E}_{k_j}|u| \right)(x) \leq C, \qquad x \in \Omega, $$
and the result follows. \finpreuve

\medskip

The next Lemma contains half of Proposition $
\refe{lemmeCalderonZygmund} $.

\begin{lemm} \label{fonctionsCalderonZygmund}  For all $ u \in L^1 (\Omega,
d\nu) $ and all $ \lambda > 0 $, we can find $ n_0 \in \Na $, an
admissible family $ {\mathcal P}  $  of partitions of $ \Omega $,
a set $ {\mathcal I} \subset \{ (k,l) \in \Za \times \Na \ | \ k
\geq - n_0 \} $ and functions $ (w_{k,l})_{(k,l) \in {\mathcal I}}
$ and $ v $ satisfying
\begin{eqnarray}
u = v + \sum_{{\mathcal I}} w_{k,l} , \qquad \qquad \qquad
\label{decou} \\
|v(x)| \leq C_3 \lambda, \qquad \mbox{a.e.} , \qquad \qquad \label{linfini} \\
\int_{\Omega} w_{k,l} \ d\nu = 0  \ \ \mbox{and} \ \
\emph{supp} \ w_{k,l} \in {\mathcal P}_{l}(k) ,  \label{intsuppw} \\
 \sum_{{\mathcal I}} \nu ({\mathcal P}_{l}(k)) \leq
 \lambda^{-1}||u||_{L^1(\Omega,d \nu)}, \qquad \label{majomesure} \\
 ||v||_{L^1(\Omega,d\nu)} + \sum_{ {\mathcal I} }
||w_{k,l}||_{L^1(\Omega,d\nu)} \leq 3
||u||_{L^1(\Omega,d\nu)} . \label{majonorml1}
\end{eqnarray}
The constant $ C_3 $ in  \reff{linfini}  is the one chosen in
\reff{constante2bis}.
\end{lemm}

\noindent {\it Proof.}  We first choose $ n_0 \in \Na $ such that
$ C_2^{1-n} 2^{n_0(n-1)}
> \lambda^{-1} ||u||_{L^1(\Omega,d\nu)} $, using the same constant $C_2$ as in
\reff{minoration1}, and then consider the admissible family of
partitions $ {\mathcal P} = \left( ({\mathcal P}_{l}(k))_{l \in
\Na} \right)_{k \geq - n_0} $ defined by  \reff{exemple2}. By
\reff{minoration1}, we have
\begin{eqnarray}
 \nu ( {\mathcal P}_{l}(-n_0) ) > \lambda^{-1}
||u||_{L^1(\Omega,d\nu)} , \label{initialise}
\end{eqnarray}
 for all $ l \in \Na $.
Next, we define $ I_{1-n_0} \subset \Na $ as the subset of
indices $ l $ such that
$$  \frac{1}{\nu ( {\mathcal P}_{l}(1-n_0) ) }
\int_{{\mathcal P}_{l}(1-n_0)}|u| \ d \nu \geq \lambda
,
 $$
and we also set $ B_{1-n_0} = \sqcup_{ l \in I_{1-n_0}}
{\mathcal P}_{l}(1-n_0) $. By induction, we then construct $
I_{k} $ and $ B_k $, for $ k \geq 2 - n_0 $, by
\begin{eqnarray}
 I_k = \left\{ l \in \Na \ \big|  \  \int_{{\mathcal
P}_{l}(k)}|u| \ d \nu \geq \lambda
 \nu ( {\mathcal P}_{l}(k) ) \ \ \mbox{and} \ \ {\mathcal P}_{l}(k) \cap B_{k-1} = \emptyset
  \right\} \nonumber \\
B_k =  B_{k-1}\sqcup \bigsqcup_{l \in I_k} {\mathcal
P}_{l}(k) . \qquad \qquad \qquad \qquad \qquad
 \nonumber
\end{eqnarray}
We now define $ B= \cup_{k \geq 1 - n_0} B_k $, $ A = \Omega
\setminus B $ and  $ {\mathcal I} = \cup_{k \geq 1 - n_0} \{ k \}
\times I_k $. Then, we first observe that for all $  k \geq 1 -
n_0 $ and all $ l \in \Na $
\begin{eqnarray}
{\mathcal P}_{l}(k) \cap A \ne \emptyset \ \ \Rightarrow \ \
\frac{1}{\nu({\mathcal P}_{l}(k))} \int_{{\mathcal
P}_{l}(k)} |u| \ d\nu < \lambda . \label{impl}
\end{eqnarray}
 Indeed, assume that $ {\mathcal P}_{l}(k) \cap A \ne \emptyset $. Then
  $ {\mathcal P}_{l}(k) \cap B_k = \emptyset $, otherwise
we could find $ k^{\prime} \leq k  $ and $ l^{\prime} \in
I_{k^{\prime}} $ such that $ {\mathcal P}_{l}(k) \cap {\mathcal
P}_{l^{\prime}}(k^{\prime}) \ne \emptyset $ in which case we would
have $ {\mathcal P}_{l}(k) \subset {\mathcal
P}_{l^{\prime}}(k^{\prime}) $ (since $ {\mathcal P} $ is non
decreasing) and thus $ {\mathcal P}_{l}(k) \cap A = \emptyset $.
Moreover,  $ l \notin I_k $ otherwise $ {\mathcal P}_{l}(k)
\subset B_k \subset B $. In addition, if $ k \geq 2 - n_0 $, $
{\mathcal P}_{l}(k) \cap B_{k-1} = \emptyset $ (since $ B_{k-1}
\subset B_k $) so the right hand side of  \reff{impl}  holds.
Then, by Lemma  \refe{standard},  \reff{impl}  shows that $ |u|
\leq \lambda $ almost everywhere on $ A $ and we set
\begin{eqnarray}
 v (x) = u (x), \qquad x \in A . \nonumber
\end{eqnarray}
  We then define $ v $ on $ B $
by
$$ v (x)  = \sum_{(k,l) \in {\mathcal I}}   \frac{1}{\nu ( {\mathcal P}_{l}(k) ) }
\int_{{\mathcal P}_{l}(k)}u \ d \nu \ \chi_{{\mathcal
P}_{l}(k)}(x), \qquad x \in B ,
$$
and, for each $ (k,l) \in {\mathcal I} $, we  define $ w_{k,l} $
on $ \Omega $ by
$$ w_{k,l} = (u-v) \chi_{{\mathcal P}_{l}(k)} . $$
With such a choice, \reff{intsuppw} holds. Observing that, by
construction, we have $ B = \sqcup_{(k,l) \in {\mathcal I}}
{\mathcal P}_{l}(k)$, one easily sees that \reff{decou} and
\reff{majonorml1} hold, using for the latter the fact that
$$ \int_{\Omega} |v| \ d \nu \leq \int_{\Omega} |u| \ d\nu . $$
 Furthermore, for all $
(k,l) \in {\mathcal I} $, we have
$$ \lambda \nu({\mathcal P}_{l}(k)) \leq \int_{{\mathcal
P}_{l}(k)} |u| \ d\nu  ,$$ and summing these estimates yields
\reff{majomesure}. We still have to check \reff{linfini}.  Since
we already know that $ |v|=|u| \leq \lambda $ a.e. on $ A $, it
only remains to show that $ |v| \leq C_3 \lambda  $ a.e. on $ B $.
Using \reff{constante2bis}, it is enough to show that, for all $
(k,l) \in {\mathcal I} $,
\begin{eqnarray}
  \int_{{\mathcal
P}_{l^{\prime}}(k-1)} |u| \ d \nu <  \lambda \nu({\mathcal
P}_{l^{\prime}}(k-1)) . \label{majointe}
\end{eqnarray}
 If $ k = 1 - n_0 $, this
follows from  \reff{initialise}. If $ k \geq 2 - n_0 $, we remark
that $ l^{\prime} \notin I_{k-1} $ otherwise $ l $ could not
belong to $ I_k $ since we would have $ {\mathcal P}_{l}(k)
\subset {\mathcal P}_{l^{\prime}}(k-1) \subset B_{k-1} $. Hence,
either \reff{majointe} holds or $ {\mathcal P}_{l^{\prime}}(k-1)
\cap B_{k-2} \ne \emptyset $, but latter is excluded, otherwise
Remark $ 3 $ (before Definition \refe{definitionremarque}) would
imply that $ {\mathcal P}_l (k) \subset {\mathcal
P}_{l^{\prime}}(k-1) \subset B_{k-2} \subset B_{k-1} $ which would
prevent $ l $ to belong to $ I_k $. The proof is complete.
\finpreuve

\bigskip

\noindent {\bf Proof of Proposition 
\ref{lemmeCalderonZygmund}.} Relabelling   the family $
({\mathcal
  P}_l(k))_{(k  , l) \in {\mathcal I}}  $  obtained in  Lemma
 \refe{fonctionsCalderonZygmund}   as $ (Q_j)_{j \in \Na} $ and
defining correspondingly the functions $ \tilde{u}  =v  $ and $
u_j = w_{k,l}  $, Lemma  \refe{fonctionsCalderonZygmund} states
that $ \reff{sommeuj} $,  $ \reff{borneponctuelle} $, $
\reff{moyennenulle}  $, $ \reff{mesureL1}  $ and $ \reff{borneL1}
$ hold. Let $ j \in \Na $. Then, there exist $ k \geq - n_0 $ and
$ l \in \Na $ such that $ Q_j = {\mathcal P}_l (k) $ (which is of
the form  \reff{exemple2} for some  $ i \in \Na $ and $ m \in
\Za^{n-1} $) and hence
$$ Q_j \subset \left\{ |r-2^{-k_+}i| \leq 2^{-k_+} \ \mbox{and} \ \frac{|\theta-2^{-k}
w([2^{-k_+}i])m|}{w(2^{-k_+ i})} \leq 2^{-k} (n-1)^{1/2}
\frac{w([2^{-k_+}i])}{w(2^{-k_+}i)} \right\} .
$$
Here the factor $ (n-1)^{1/2} $ is due to the inclusion of all
cubes of side $2$ into euclidean balls of radius $ (n-1)^{1/2} $.
By setting
$$ r_j = 2^{-k_+}i , \qquad \theta_j = 2^{-k} w([2^{-k_+}i]) m $$
and using the fact that $  w([2^{-k_+}i])/w(2^{-k_+}i) \leq C $,
which is due to $ \reff{diag} $, we see that the first inclusion
in $ \reff{tjun} $ holds with
$$ t_j = 2^{- k_+} + C 2^{-k}(n-1)^{1/2} $$
if this quantity is $ \leq 1 $, and otherwise, if it is $ > 1 $,
that the first inclusion of $ \reff{tjde} $ holds with
$$ t_j = \max ( 2 , 2^{- k_+} + C 2^{-k}(n-1)^{1/2} ) . $$
Note that $ t_j \approx 2^{-k} $ in all cases. Fix now $ D > 1 $.
If $ t_j \leq 1 $, we define $ Q_j^* $ as the third set in
\reff{tjun}, and if $ t_j > 1 $ (ie $ t_j \geq 2 $), we define $
Q_j^* $ as the third set of  \reff{tjde} so that, in both cases,
the last inclusions of \reff{tjun} and \reff{tjde} are actually
equalities. It thus remains to prove \reff{tjunb}. If $ t_j  \leq
1  $, then, using  \reff{diagequivalent}, we have
$$ \nu (Q_j^*) \lesssim \int_{r_j - D t_j}^{r_j + D t_j} \left(D t_j \frac{w(r_j)}{w(r)} \right)^{n-1}
d r \lesssim  C^{n-1} e^{(n-1)CD} (D t_j)^{n} \approx C_D 2^{-kn}
\approx C_D  \nu (Q_j) , $$ using \reff{minoration2} for the last
estimate. Similarly, if  $ t_j  > 1 $, then
$$   \nu (Q_j^*) \lesssim \int_{r_j - 2}^{r_j + 2} \left(D t_j \frac{w(r_j)}{w(r)} \right)^{n-1}
d r \lesssim 4  C^{n-1} e^{2(n-1)C} (D t_j)^{n-1} \approx C_D
2^{-(n-1)k} \approx C_D \nu (Q_j) ,  $$ using,
 for the last estimate, \reff{minoration1} if $ k \leq 0 $ or \reff{minoration2} for the finite number of
 $ k \geq 1 $ such that $ 2^{-k_+} + C 2^{-k}(n-1)^{1/2} > 1 $. The proof is complete. \finpreuve


\begin{thebibliography}{99}

\bibitem{Bo2} {\sc J. M. Bouclet} {\it Semi-classical functional calculus on 
    manifolds with ends and weighted $L^p$ estimates}, preprint.

\bibitem{Bohyperbolique} {\sc \name}, {\it Strichartz estimates on asymptotically hyperbolic manifolds}, preprint.

\bibitem{BoTz1} {\sc J. M. Bouclet, N. Tzvetkov},
{\it Strichartz estimates for long range perturbations}, Amer. J.
Math., to appear.

\bibitem{BoTz2} {\sc \name}, {\it On global Strichartz estimates for non
    trapping metrics}, J. Funct. Analysis,  to appear.

\bibitem{BGT1} {\sc N. Burq, P. G\'erard, N. Tzvetkov},
{\it Strichartz inequalities and the non linear Schr\"odinger
equation on compact manifolds},  Amer. J. Math. 126, 569-605
(2004).

\bibitem{Coul1} {\sc T. Coulhon, X. T. Duong, X. D. Li}, {\it
    Littlewood-Paley-Stein functions on complete Riemannian manifolds for $ 1
    \leq p \leq 2 $}, Studia Math. 154, no. 1, 37-57 (2003).







\bibitem{KlRo1} {\sc S. Klainerman, I. Rodnianski}, {\it A geometric approach to the
Littlewood-Paley theory}, Geom. Funct. Anal.  16 No. 1, 126-163 (2006).

\bibitem{Loho1} {\sc N. Lohou\'e}, {\it Estimations des fonctions de
    Littlewood-Paley-Stein sur les vari\'et\'es \`a courbure non positive},
  Ann. Sci. \'Ecole Norm. Sup. 20, 505-544 (1987).




\bibitem{OlZh1} {\sc G. Olafsson, S. Zheng}, {\it Harmonic analysis related to
  Schr\"odinger operators}, arXiv:0711.3262v1.

\bibitem{Schl1} {\sc W. Schlag}, {\it A remark on Littlewood-Paley theory for
    the distorted Fourier transform}, Proc. Amer. Math. Soc. 135 no. 2,
  437-451 (2007).

\bibitem{Sogg1} {\sc C. D. Sogge},
{\it Fourier integrals in classical analysis}, Cambridge tract in
Mathematics, Cambridge Univ. Press (1993).

\bibitem{Stei1} {\sc E. M. Stein},
{\it Singular integrals and differentiability properties of
functions}, Princeton Univ. Press (1970).

\bibitem{Stei2} {\sc \name}, {\it Topics in harmonic analysis related to the
    Littlewood-Paley theory}, Princeton Univ. Press (1970).




\bibitem{Tayl0} {\sc M. Taylor}, {\it $L^p$ estimates on functions of the Laplace
operator}, Duke Math. J. Vol. 58, No. 3, 773-793 (1989).


\bibitem{Tayl1} {\sc \name},
{\it  Partial Differential Equations III, Nonlinear Equations},
Appl. Math. Sci. 117, Springer, (1996).

\end{thebibliography}
\end{document}